\documentclass[lettersize,journal]{IEEEtran}
\usepackage{amsmath,amsfonts}
\usepackage{amssymb}  % assumes amsmath package installed
\usepackage{algorithm}
\usepackage{array}
\usepackage[caption=false,font=normalsize,labelfont=sf,textfont=sf]{subfig}
\usepackage{textcomp}
\usepackage{stfloats}
\usepackage{url}
\usepackage{verbatim}
\usepackage{graphicx}
\usepackage{cite}
\usepackage{enumerate}
\usepackage{enumitem}
\usepackage{nomencl}
\makenomenclature
\usepackage{etoolbox}

\newcounter{Ax}

\usepackage{color}

\usepackage{caption}
\usepackage{multicol}
\usepackage{multirow}
\usepackage{diagbox}
\usepackage{makecell}

\usepackage{booktabs}
\usepackage[noend]{algpseudocode}
\usepackage{bbm}
\usepackage{dsfont}

\newtheorem{lemma}{Lemma}
\newtheorem{assumption}{Assumption}
\newtheorem{remark}{Remark}
\newtheorem{theorem}{Theorem}

\allowdisplaybreaks[4]  % Allow the equations to continue in a new page

\begin{document}

\title{A Hierarchical OPF Algorithm with Improved Gradient Evaluation in Three-Phase Networks}

\author{Heng Liang, Xinyang Zhou,~\IEEEmembership{Member,~IEEE}, and Changhong Zhao,~\IEEEmembership{Senior Member,~IEEE}
        % <-this % stops a space

\thanks{This work was supported by the Hong Kong Research Grants Council through ECS Award No. 24210220. (Corresponding author: Changhong Zhao.)}
\thanks{H. Liang and C. Zhao are with the Department of Information Engineering, the Chinese University of Hong Kong, New Territories, Hong Kong SAR, China
(e-mail:  lh021@ie.cuhk.edu.hk; chzhao@ie.cuhk.edu.hk).}% <-this % stops a space
\thanks{X. Zhou is with the National Renewable Energy Laboratory, Golden, CO, USA
(e-mail: xinyang.zhou@nrel.gov).}
\thanks{This work has been submitted to the IEEE for possible publication. Copyright may be transferred without notice, after which this version may no longer be accessible.}}



\maketitle

\begin{abstract}

Linear approximation commonly used in solving alternating-current optimal power flow (AC-OPF) simplifies the system models but incurs accumulated voltage errors in large power networks. Such errors will make the primal-dual type gradient algorithms converge to the solutions at which the power networks may be exposed to the risk of voltage violation. In this paper, we improve a recent hierarchical OPF algorithm that rested on primal-dual gradients evaluated with a linearized distribution power flow model. Specifically, we propose a more accurate gradient evaluation method based on a three-phase unbalanced nonlinear distribution power flow model to mitigate the errors arising from model linearization. The resultant gradients feature a blocked structure that enables us to further develop an improved hierarchical primal-dual algorithm to solve the OPF problem. Numerical results on the IEEE $123$-bus test feeder and a $4,518$-node test feeder show that the proposed method can enhance the overall voltage safety while achieving comparable computational efficiency with the linearized algorithm.

\end{abstract}


\printnomenclature 
\nomenclature[01]{$\mathcal{N}$}{Set of buses.}
\nomenclature[02]{$\mathcal{N}^{+}$}{$\mathcal{N}^{+}:=\mathcal{N}\cup\{0\}$, set of buses including the slack bus.}
\nomenclature[03]{$\mathcal{E}$}{Set of lines.}
\nomenclature[04]{$\Phi_j$}{Phase set of bus $j$.}
\nomenclature[05]{$\Phi_{ij}$}{Phase set of line $(i,j)$.}
\nomenclature[06]{$V_j^{\phi}$}{Complex phase voltage of bus $j$ at phase $\phi$.}
\nomenclature[07]{$V_{j}$}{$V_{j}:=(V_{j}^{\phi},\phi\in\Phi_{j})$, column vector of complex voltage of bus $j$ at all phases.}
\nomenclature[08]{$v_{j}$}{$v_{j}:=V_{j}V_{j}^{H}$, matrix of complex squared voltage.}
\nomenclature[09]{$I_{ij}^{\phi}$}{Complex current of line $(i,j)$ at phase $\Phi_{ij}$.}
\nomenclature[10]{$I_{ij}$}{$I_{ij}:=(I_{ij}^{\phi},\phi\in\Phi_{ij})$, column vector of complex current of line $(i,j)$ at all phases.}
\nomenclature[11]{$\ell_{ij}$}{$\ell_{ij}:=I_{ij}I_{ij}^{H}$, matrix of complex squared current.}
\nomenclature[12]{$S_{ij}$}{$S_{ij}:=V_{i}^{\Phi_{ij}}I_{ij}^{H}$, matrix of power flow at the sending-end of line $(i,j)$.}
\nomenclature[13]{$\Lambda_{ij}$}{Column vector of sending-end three-phase power flow.}
\nomenclature[16]{$z_{ij}$}{Matrix of series impedance of line $(i,j)$.}
\nomenclature[17]{$p_j^{\phi}/q_j^{\phi}$}{Active/Reactive power injection of bus $j$ at phase $\phi$.}
\nomenclature[18]{$s_j^{\phi}$}{Complex power injection of bus $j$ at phase $\phi$.}
\nomenclature[19]{$s_j$}{$s_{j}:=(s_{j}^{\phi},\phi\in\Phi_{j})$, column vector of complex power injection of bus $j$.}
\nomenclature[20]{$u_j^{\phi}$}{$u_j^{\phi}:=[p_j^{\phi},q_j^{\phi}]^{\top}$, power injection of bus $j$ at phase $\phi$.}
\nomenclature[21]{$\boldsymbol{p}$}{$\boldsymbol{p}:=[[p_1^{\phi}]^{\top}_{\phi\in\Phi_{1}},\dots,[p_N^{\phi}]^{\top}_{\phi\in\Phi_{N}}]^{\top}$, three-phase active power injection vector.}
\nomenclature[22]{$\boldsymbol{q}$}{$\boldsymbol{q}:=[[q_1^{\phi}]^{\top}_{\phi\in\Phi_{1}},\dots,[q_N^{\phi}]^{\top}_{\phi\in\Phi_{N}}]^{\top}$, three-phase reactive power injection vector.}
\nomenclature[23]{$\boldsymbol{v}$}{$\boldsymbol{v}:=[[v_1^{\phi\phi}]^{\top}_{\phi\in\Phi_{1}},\dots,[v_N^{\phi\phi}]^{\top}_{\phi\in\Phi_{N}}]^{\top}$, three-phase squared voltage magnitude vector.}
\nomenclature[24]{$\underline{\boldsymbol{v}}/\overline{\boldsymbol{v}}$}{Lower/Upper bounds of three-phase squared voltage magnitude.}
\nomenclature[25]{$\underline{\boldsymbol{\mu}}/\overline{\boldsymbol{\mu}}$}{Dual variables for the lower/upper bounds of three-phase squared voltage magnitude.}
\nomenclature[26]{$\gamma$}{Constant matrix for balanced voltage approximation.}
\nomenclature[50]{$\overline{x}$}{Conjugate of complex number $x$.}
\nomenclature[51]{$\hat{v}$}{Voltages of linear approximation of BFM (other variables in BFM are similar, and the below is the same).}
\nomenclature[52]{$\widetilde{v}$}{Voltages of improved gradient evaluation of BFM.}
\nomenclature[53]{$v^{\phi\psi}_i$}{Element of complex matrix $v_i$ at row $\phi$, column $\psi$.}

\section{Introduction}
\label{sec::introduction}

\IEEEPARstart{O}{ptimal} power flow (OPF) is a fundamental optimization problem that aims to find a cost-minimizing operating point subject to the constraints of physical laws and safety limits. OPF underlies many important power system applications such as demand response, state estimation, unit commitment, and voltage regulation. With the increasing penetration of controllable units (smart appliances, electric vehicles, energy storage devices, etc.) and thus the growing size of OPF problems, solving OPF is becoming more challenging with heavy computations and intense communications. Especially, the growth of wind and solar generations introduces increased variations to power systems and calls for timely response and real-time optimization. This trend poses a rising need for fast and scalable OPF solvers, which is especially urgent in power distribution networks where massive renewable energy resources and controllable units are being deployed. However, distribution networks are also where the algorithm speed and scalability requirements are most difficult to meet, as the high resistance-to-reactance ratios of distribution lines necessitate the usage of nonlinear and nonconvex alternate-current (AC) power flow models rather than their simple direct-current approximations.

Numerous efforts have been made to overcome this challenge. Many of them conduct convex relaxation, such as semidefinite program (SDP) relaxation \cite{bai2008semidefinite} and second order cone program (SOCP) relaxation \cite{farivar2013branch}; see \cite{low2014convex,low2014convexii} for their connection and equivalence. Meanwhile, convex inner approximation and linearization of AC power flow were proposed to simplify the system models. The classical linearized distribution power flow model was derived by neglecting the non-linear line loss terms in \cite{baran1989optimalC,baran1989optimalCii}, while linearization around a fixed-point was shown to have a better approximation accuracy \cite{Bernstein2018loadflow,ESchweitzer2020lossy,RCheng2022onlinelinear}. A comprehensive review of approximate models for AC power flow was provided in \cite{molzahn2019survey}. Based on these models, the existing methods often rely on off-the-shelf solvers to solve OPF, which often suffer from overwhelming computations and communications in large-scale networks. Distributed OPF algorithms were designed and shown to be more scalable in terms of computation and more robust to single point failure, compared to their centralized counterparts \cite{dall2013distributed,erseghe2014distributed,zhang2014optimal}. The distributed algorithms are typically implemented in two ways. The first way relies on information exchange between the neighboring buses \cite{peng2016distributed,peng2015distributed,SMagnusson2020distributed,NPatari2022distributed}. The second way involves a central controller (CC), with all buses sending/receiving information to/from the CC \cite{tang2017distributed,dall2016optimal,zhou2019hierarchical,zhou2019accelerated}. Due to the complex interdependence and nonconvexity of the nonlinear AC power flow, most distributed OPF algorithms utilized linear approximation to facilitate the design process. To further reduce computational efforts associated with solving AC power flow, some OPF algorithms were implemented by iteratively actuating the power system with intermediate decisions and updating the decisions based on system feedback \cite{gan2016online,bolognani2014distributed,bernstein2019real}.

From the vast literature, we bring attention to a hierarchical distributed primal-dual gradient algorithm \cite{zhou2019hierarchical} and its extension to three-phase unbalanced networks \cite{zhou2019accelerated}. This algorithm leveraged the radial structure of distribution networks to avoid repetitive computation and communication, and thus significantly accelerated the solution process of large-scale OPF problems. However, the gradient used in this algorithm was derived from the linearized distribution power flow model \cite{baran1989optimalC}. Such linearization, as we will analyze, may cause the solver to optimistically estimate nodal voltages to be safe, while they actually already violate safety limits.

To prevent such violation, we develop an improved gradient evaluation method inspired by the chain rule of derivatives. Take the nonlinear single-phase branch flow model as an example \cite{baran1989optimalC,farivar2013branch}. The variables in that model can be grouped into three sets, which are the active and reactive power injections $\boldsymbol{u}:=(\boldsymbol{p},\boldsymbol{q})$, the squared voltage magnitudes $\boldsymbol{v}(\boldsymbol{u})$ and the branch power flows and squared currents $(\boldsymbol{P},\boldsymbol{Q},\boldsymbol{\ell})$. Then the partial derivative of squared voltage magnitudes over power injections can be taken as:
\begin{alignat}{2}
\frac{\partial \boldsymbol{v}(\boldsymbol{u})}{\partial \boldsymbol{u}}=\frac{\partial \boldsymbol{v}(\boldsymbol{u})}{\partial (\boldsymbol{P},\boldsymbol{Q},\boldsymbol{\ell})}\cdot \frac{\partial (\boldsymbol{P},\boldsymbol{Q},\boldsymbol{\ell})}{\partial \boldsymbol{u}}, \nonumber
\end{alignat}
where we replace the second part with the partial derivatives $ \frac{\partial (\hat{\boldsymbol{P}},\hat{\boldsymbol{Q}},\hat{\boldsymbol{\ell}})}{\partial \boldsymbol{u}}$ calculated from the linearized model $(\hat{\boldsymbol{P}}(\boldsymbol{u}),\hat{\boldsymbol{Q}}(\boldsymbol{u}),\hat{\boldsymbol{\ell}}(\boldsymbol{u}))$. In particular, the proposed method preserves an approximate calculation associated with the quadratic terms in the nonlinear AC power flow model, instead of disregarding them. It is worth pointing out that this method is not only applicable to single-phase networks, as shown in our preliminary conference paper \cite{liang2022hierarchical}, but can also be extended to three-phase unbalanced networks as introduced in this paper. Our analysis shows that with moderate extra computations, the proposed method returns more accurate gradient evaluations than \cite{zhou2019hierarchical,zhou2019accelerated}, while preserving the blocked gradient structure therein that enables us to develop a scalable and voltage-safe hierarchical OPF algorithm. To summarize, the contributions of this paper are:
\begin{itemize}
    \item We propose an improved gradient evaluation method for three-phase unbalanced nonlinear AC power flow, which extends our previous work \cite{liang2022hierarchical} on a single-phase model.
    \item Based on the improved gradient evaluation, we design a voltage-safe hierarchical OPF algorithm for three-phase unbalanced power distribution networks.
    \item We perform a rigorous voltage error analysis between the nonlinear distribution power flow model and its linear approximation. A convergence analysis of the proposed algorithm is presented to show the smaller sub-optimality gap of the result with more accurate gradient evaluation.
    \item Numerical experiments on large-scale three-phase unbalanced networks, including the IEEE $123$-bus test feeder and a $4,518$-node test feeder, are performed to demonstrate the enhanced voltage safety and computational efficiency of the proposed method.
\end{itemize}

The rest of this paper is organized as follows. Section \ref{sec:multi:model} introduces a three-phase unbalanced distribution network model, an OPF, and a primal-dual algorithm to solve it. Section \ref{sec::improved} motivates and proposes the improved gradient evaluation method. Section \ref{sec:multi:algorithm} designs an improved hierarchical OPF algorithm and analyzes its convergence. Section \ref{sec:numerical} reports the numerical experiments. Section \ref{sec:conclusion} concludes this paper.

\section{Modeling and Preliminary Algorithm}
\label{sec:multi:model}

% This section defines a three-phase OPF problem and proposes an improved gradient evaluation method for three-phase unbalanced networks.

\subsection{Branch Flow Model and OPF Formulation}

We denote the set of complex numbers, the set of $n$-dimensional complex vectors, and the set of $m\times n$ complex matrices by $\mathbb{C}$, $\mathbb{C}^{n}$, $\mathbb{C}^{m\times n}$, respectively. Let $\overline{x}$ denote the conjugate of a complex number $x\in \mathbb{C}$. Let $\left|\cdot\right|$ take the dimension of a vector, and $\left(\cdot\right)^{H}$ be the conjugate transpose of a matrix or vector.

We model a three-phase unbalanced distribution power network as a directed tree graph $\mathcal{T}:=\{\mathcal{N}^{+}, \mathcal{E}\}$, where $\mathcal{N}^{+}=\{0\}\cup \mathcal{N}$, with bus $0$ indexing the slack (root) bus and $\mathcal{N}=\{1,...,N\}$ containing other buses. All the lines $(i,j) \in \mathcal{E}$ have reference directions that point from the root towards leaves, so that each bus $j\in \mathcal{N}$ has a unique upstream bus $i$, while the root bus has no upstream bus. Let $a$, $b$, $c$ denote the three phases. We use $\Phi_{j}$ to denote the phase set of bus $j$ and $\Phi_{i j}$ the phase set of line $(i,j)$, e.g., $\Phi_{j}=\{a,b,c\}$ if bus $j$ has all the three phases, and $\Phi_{ij}=\{a\}$ if line $(i,j)$ has phase ``$a$'' only. 

For each phase $\phi\in \Phi_{j}$ at bus $j \in \mathcal{N}^{+}$, let $V_{j}^{\phi}\in \mathbb{C}$ be the complex voltage, and $s_{j}^{\phi}:=p_{j}^{\phi}+\mathbf{i}q_{j}^{\phi}$ be the complex power injection (i.e., power generation minus consumption at buses $j \in\mathcal{N}$ and power flow from the upper grid into the distribution network at the slack bus). 
Define $V_{j}:=(V_{j}^{\phi},\phi\in\Phi_{j}) \in \mathbb{C}^{\left|\Phi_{j}\right|}$, $s_{j}:=(s_{j}^{\phi},\phi\in\Phi_{j})\in \mathbb{C}^{\left|\Phi_{j}\right|}$ and $v_{j}:=V_{j}V_{j}^{H}\in \mathbb{C}^{\left|\Phi_{j}\right|\times \left|\Phi_{j}\right|}$. 

For each line $(i,j)\in \mathcal{E}$, let $I_{ij}^{\phi} \in \mathbb{C}$ be the complex current on phase $\phi \in \Phi_{i j}$, and $z_{ij} \in \mathbb{C}^{\left|\Phi_{ij}\right| \times\left|\Phi_{ij}\right|}$ be the series impedance matrix. Define $I_{ij}:=(I_{ij}^{\phi},\phi\in\Phi_{ij}) \in \mathbb{C}^{\left|\Phi_{ij}\right|}$, $S_{ij}:=V_{i}^{\Phi_{ij}}I_{ij}^{H} \in \mathbb{C}^{\left|\Phi_{ij}\right|\times \left|\Phi_{ij}\right|}$ (where $V_{i}^{\Phi_{ij}}$ denotes the subvector of $V_i$ on phases $\Phi_{ij}$), and $\ell_{ij}:=I_{ij}I_{ij}^{H} \in \mathbb{C}^{\left|\Phi_{ij}\right|\times \left|\Phi_{ij}\right|}$. 

Let $\operatorname{diag}(\cdot)$ denote the column vector composed of diagonal elements of a matrix; in the other way, 
$\operatorname{Diag}(\cdot)$ converts a column vector into a diagonal matrix.
Let $v^{\phi\psi}_{i}$ for $\phi,\psi \in \Phi_{i}$ denote the element of complex matrix $v_{i}$ at row $\phi$, column $\psi$. Phases $a$, $b$, $c$ are represented by numbers 0, 1, 2, respectively, wherever needed.

Consider the nonlinear three-phase unbalanced distribution branch flow model (BFM) \cite{gan2014convex,zhao2017optimal}:
\begin{subequations}\label{distflow} 
\begin{alignat}{2} 
&v_{j}= v_{i}^{\Phi_{ij}}-\left(S_{i j} z_{i j}^{H}+z_{i j} S_{i j}^{H}\right)+z_{i j} \ell_{i j} z_{i j}^{H},  \forall (i, j)\in \mathcal{E}, \label{disflow::v}\\
& \operatorname{diag}\left(S_{i j}-z_{i j} \ell_{i j}\right)-\sum_{k:(j,k)\in\mathcal{E}} \operatorname{diag}\left(S_{j k}\right)^{\Phi_j}=-s_{j}, \nonumber\\
& \qquad\qquad \qquad\qquad\qquad\qquad\qquad\qquad\qquad ~\forall j \in \mathcal{N}^+,\label{disflow::PQ}\\
&\left[\begin{array}{cc}
v_{i}^{\Phi_{ij}} & S_{i j} \\
S_{i j}^{H} & \ell_{i j}
\end{array}\right] \succeq 0, \quad \forall (i,j) \in \mathcal{E}, \label{disflow::PSD}\\
&\operatorname{rank}\left(\left[\begin{array}{cc}
v_{i}^{\Phi_{ij}} & S_{i j} \\
S_{i j}^{H} & \ell_{i j}
\end{array}\right]\right)=1, \quad \forall (i, j) \in \mathcal{E}. \label{disflow::rank}
\end{alignat}
\end{subequations}

Suppose voltage $V_0$ and $v_0$ at the slack bus are given and fixed. We use $[p_i^{\phi}]_{\phi \in \Phi_i}$ to denote the column vector of active power injections across all the phases of bus $i$, and so on. Define $\boldsymbol{p}:=[[p_1^{\phi}]^{\top}_{\phi\in\Phi_{1}},\dots,[p_N^{\phi}]^{\top}_{\phi\in\Phi_{N}}]^{\top}$, $\boldsymbol{q}:=[[q_1^{\phi}]^{\top}_{\phi\in\Phi_{1}},\dots,[q_N^{\phi}]^{\top}_{\phi\in\Phi_{N}}]^{\top}$ as the three-phase power injection vectors, and $\boldsymbol{v}:=[[v_1^{\phi\phi}]^{\top}_{\phi\in\Phi_{1}},\dots,[v_N^{\phi\phi}]^{\top}_{\phi\in\Phi_{N}}]^{\top}$ as the three-phase squared voltage magnitude vector. The function $\boldsymbol{v}(\boldsymbol{p},\boldsymbol{q})$ is implicitly well defined by the nonlinear BFM (\ref{distflow}) under normal operating conditions \cite{cwang2018explicit}. Consider the following three-phase OPF problem:
\begin{subequations} \label{MultiOPF}
\begin{alignat}{2}
\min _{\boldsymbol{p}, \boldsymbol{q}} & \sum_{i \in \mathcal{N}}\sum_{\phi \in \Phi_{i} } f_{i}^{\phi}\left(p_{i}^{\phi}, q_{i}^{\phi}\right) \label{MultiOPF:obj}\\
\text { s.t. } & \underline{\boldsymbol{v}} \leqslant \boldsymbol{v}(\boldsymbol{p}, \boldsymbol{q}) \leqslant \overline{\boldsymbol{v}},\label{MultiOPF:v}\\
& \left(p_{i}^{\phi}, q_{i}^{\phi}\right) \in \mathcal{Y}_{i}^{\phi}, \quad\forall \phi \in \Phi_{i},\forall i \in \mathcal{N},
\end{alignat}
\end{subequations}
where $f_{i}^{\phi}\left(p_{i}^{\phi}, q_{i}^{\phi}\right)$ is a strongly convex cost function (e.g., a quadratic function) of the controllable power injection at phase $\phi$ of bus $i$. The power injections are confined by \textit{compact convex} set $\mathcal{Y}_{i}^{\phi}$ for each $i\in \mathcal{N}$, $\phi\in \Phi_i$, for instance, a box:
\begin{eqnarray} \nonumber
\mathcal{Y}_{i}^{\phi}=\left\{\left(p_{i}^{\phi}, q_{i}^{\phi}\right) \mid \underline{p}_{i}^{\phi} \leqslant p_{i}^{\phi} \leqslant \bar{p}_{i}^{\phi}, \underline{q}_{i}^{\phi} \leqslant q_{i}^{\phi} \leqslant \bar{q}_{i}^{\phi}\right\}.
\end{eqnarray}

\subsection{Primal-dual Gradient Algorithm}
Let $\underline{\boldsymbol{\mu}}$ and $\overline{\boldsymbol{\mu}}$ be the dual variables associated with the left-hand-side (LHS) and right-hand-side (RHS) of (\ref{MultiOPF:v}), respectively. To design a convergent primal-dual algorithm, we consider the regularized Lagrangian of OPF problem (\ref{MultiOPF}):
\begin{alignat}{2}\label{MultiRLagarangian}  \nonumber
\mathcal{L}_{\epsilon}(\boldsymbol{p}, \boldsymbol{q} ; \overline{\boldsymbol{\mu}}, \underline{\boldsymbol{\mu}})= \sum_{i \in \mathcal{N}}\sum_{\phi \in \Phi_{i} } f_{i}^{\phi}\left(p_{i}^{\phi}, q_{i}^{\phi}\right)-\frac{\epsilon}{2}\|\boldsymbol{\mu}\|_{2}^{2}\\ 
+\underline{\boldsymbol{\mu}}^{\top}(\underline{\boldsymbol{v}}-\boldsymbol{v}(\boldsymbol{p}, \boldsymbol{q}))+\overline{\boldsymbol{\mu}}^{\top}(\boldsymbol{v}(\boldsymbol{p}, \boldsymbol{q})-\overline{\boldsymbol{v}}),
\end{alignat}
where $\boldsymbol{\mu}:=[\underline{\boldsymbol{\mu}}^{\top},\overline{\boldsymbol{\mu}}^{\top}]^{\top}$ and $\epsilon >0$ is a regularization factor.

\begin{remark}
The OPF problem (\ref{MultiOPF}) is naturally nonconvex due to the rank constraint (\ref{disflow::rank}). It is known that function (\ref{MultiRLagarangian}) is strongly concave in dual variables $\boldsymbol{\mu}$ \cite{koshal2011multiuser}. A saddle point of (\ref{MultiRLagarangian}) serves as an approximate sub-optimal solution to OPF problem (\ref{MultiOPF}), with its error bounded in terms of $\epsilon$ \cite{tang2018feedback}. 
\end{remark} 

A primal-dual gradient algorithm to approach a saddle point of function (\ref{MultiRLagarangian}) takes the following form, for all buses $h \in\mathcal{N}$ and phases $\varphi \in \Phi_h$:
\begin{subequations} \label{MultiIter}
\begin{alignat}{2}
&\qquad \quad p_{h}^{\varphi}(t+1)= \bigg[p_{h}^{\varphi}(t)-\sigma_{u}\frac{\partial f(\boldsymbol{u}(t))}{\partial p_{h}^{\varphi}}  \nonumber\\
& \quad -\sigma_{u}\sum_{j \in \mathcal{N}} \sum_{\phi \in \Phi_{j}} \frac{\partial v_{j}^{\phi\phi} (\boldsymbol{u}(t))}{\partial p_{h}^{\varphi}} \left(\overline{\mu}_{j}^{\phi}(t)-\underline{\mu}_{j}^{\phi}(t)\right)\bigg]_{\mathcal{Y}_{h}^{\varphi}}, \label{MultiIter:p} \\
&\qquad \quad q_{h}^{\varphi}(t+1)= \bigg[q_{h}^{\varphi}(t)-\sigma_{u}\frac{\partial f(\boldsymbol{u}(t))}{\partial q_{h}^{\varphi}}  \nonumber\\
&\quad  -\sigma_{u}\sum_{j \in \mathcal{N}} \sum_{\phi \in \Phi_{j}} \frac{\partial v_{j}^{\phi\phi}(\boldsymbol{u}(t))}{\partial q_{h}^{\varphi}} \left(\overline{\mu}_{j}^{\phi}(t)-\underline{\mu}_{j}^{\phi}(t)\right)\bigg]_{\mathcal{Y}_{h}^{\varphi}}, \label{MultiIter:q} \\
&\underline{\mu}_{h}^{\varphi}(t+1)= \Big[\underline{\mu}_{h}^{\varphi}(t)+\sigma_{\mu}\big(\underline{v}_{h}^{\varphi}-v_{h}^{\varphi\varphi}(t)-\epsilon \underline{\mu}_{h}^{\varphi}(t)\big)\Big]_{+}, \label{MultiIter:umu} \\
&\overline{\mu}_{h}^{\varphi}(t+1)= \Big[\overline{\mu}_{h}^{\varphi}(t)+\sigma_{\mu}\big(v_{h}^{\varphi\varphi}(t)-\overline{v}_{h}^{\varphi}-\epsilon \overline{\mu}_{h}^{\varphi}(t)\big)\Big]_{+}, \label{MultiIter:omu}
\end{alignat}
\end{subequations}
where $\boldsymbol{u}:=[\boldsymbol{p}^{\top},\boldsymbol{q}^{\top}]^{\top}$ is the vector of controllable power injections, and $f(\boldsymbol{u}(t)):=\sum_{i\in \mathcal{N}}\sum_{\phi\in\Phi_{i} }f_{i}^{\phi}(p_{i}^{\phi},q_{i}^{\phi})$ is the objective function (\ref{MultiOPF:obj}).  The subscripts $[\cdot]_{\mathcal{Y}^{\varphi}_{h}}$ and $[\cdot]_+$ represent the projections onto the feasible power injection region $
\mathcal{Y}_{h}^{\varphi}$ and the non-negative orthant, respectively. Positive constants $\sigma_{u}$ and $\sigma_{\mu}$ are the step sizes corresponding to the primal and dual updates, respectively. 
% Note that $v_{h}^{\varphi\varphi}(t)$ can be measured from the power network after actuating the network with the intermediate decisions $\boldsymbol{u}_{\Xi}(t)$ during the process (\ref{MultiIter}). This \emph{feedback-based} implementation has been widely adopted to compensate for likely modeling errors \cite{zhou2019accelerated,gan2016online,bernstein2019real,bolognani2014distributed}.

As $\partial f/\partial p_h$ and $\partial f/\partial q_h$ are local at each bus and easy to compute, the main difficulty in implementing dynamics (\ref{MultiIter}) is the calculation of $\partial{\boldsymbol{v}}/\partial{\boldsymbol{p}}$ and $\partial{\boldsymbol{v}}/\partial{\boldsymbol{q}}$. Prior methods \cite{gan2016online,tang2017realtime} used either backward-forward sweep or matrix inverse to calculate $\partial{\boldsymbol{v}}/\partial{\boldsymbol{p}}$ and $\partial{\boldsymbol{v}}/\partial{\boldsymbol{q}}$, which involve heavy computation and may not be feasible for large-scale networks. Moreover, they are only applied to single-phase networks, not the more realistic three-phase networks. To the best of our knowledge, no existing work is using exact $\partial{\boldsymbol{v}}/\partial{\boldsymbol{p}}$ and $\partial{\boldsymbol{v}}/\partial{\boldsymbol{q}}$ based on the nonlinear three-phase BFM (\ref{distflow}). The challenges mainly lie in the complex interdependence between variables and the nonconvexity and nonlinearity induced by the rank-one constraint (\ref{disflow::rank}). A common linear approximation to the three-phase BFM, e.g., in \cite{gan2014convex,zhou2019accelerated,zhou2020multi}, assumes negligible power loss and balanced voltages to obtain a simple evaluation of $\partial{\boldsymbol{v}}/\partial{\boldsymbol{p}}$ and $\partial{\boldsymbol{v}}/\partial{\boldsymbol{q}}$. As we will analyze shortly, this linear approximation will expose the power network to the risk of voltage violation.
The voltage feedback  \cite{zhou2019accelerated,gan2016online,bernstein2019real,bolognani2014distributed} from the physical (power) network can partly resolve such violation caused by modeling error. However, ubiquitous feedback across a network is either expensive or unavailable due to limited measurement capabilities in existing systems.

\section{Improved Gradient Evaluation in Three-Phase Networks}
\label{sec::improved}
\subsection{Motivation for Improvement}
\label{subsec::motivation}

We motivate the proposed improved gradient evaluation method by introducing the limitation of an existing linearized three-phase model, in which the following approximations are made \cite{zhou2019accelerated,gan2014convex}.

     \textit{Approx. 1}: The line losses are small, i.e., $z_{ij}\ell_{ij}\ll S_{ij}$, and are thus ignored. 
    
    \textit{Approx. 2}: The voltages are nearly balanced, i.e., if $\Phi_{i}=\{a,b,c\}$, the three-phase voltage magnitudes are equal and the phase differences are close to $2\pi/3$:
    $$\frac{V_{i}^{a}}{V_{i}^{b}}\approx \frac{V_{i}^{b}}{V_{i}^{c}} \approx \frac{V_{i}^{c}}{V_{i}^{a}}\approx e^{\mathbf{i}2\pi/3}.$$

This leads to a linearized distribution network BFM:
\begin{subequations} \label{linearmodel}
\begin{alignat}{2}
&\hat{v}_{j}= \hat{v}_{i}^{\Phi_{ij}}-\left(\hat{S}_{i j} z_{i j}^{H}+z_{i j} \hat{S}_{i j}^{H}\right),  \forall (i,j)\in \mathcal{E},\\
& \hat{S}_{i j}=\gamma^{\Phi_{ij}} \operatorname{Diag}\left(\hat{\Lambda}_{i j}\right),~  \forall (i,j)\in \mathcal{E}, \label{linearmodel::S}\\
&\hat{\Lambda}_{i j}-\sum_{k: (j,k)\in \mathcal{E}}\hat{\Lambda}_{j k}^{\Phi_j}=-s_{j}, ~\forall j \in \mathcal{N}^{+},
\end{alignat}
\end{subequations}
where $\hat{\Lambda}_{i j} \in \mathbb{C}^{\left|\Phi_{ij}\right|}$ is the column vector of sending-end three-phase power flow onto the line $(i,j)\in \mathcal{E}$, and constant matrix $\gamma$ is defined as:
\begin{eqnarray}\nonumber
\gamma:=\beta \beta^{H},~\textnormal{where}~
\beta=\left[\begin{array}{c}
1 \\
\alpha \\
\alpha^{2}
\end{array}\right],\quad \alpha=e^{-\mathbf{i}2 \pi / 3}.
\end{eqnarray}
The linearized BFM (\ref{linearmodel}) yields approximate gradients:
\begin{subequations} \label{MLineargrad}
\begin{alignat}{2}
\frac{\partial \hat{\Lambda}_{i j}^{\phi}}{\partial p_{h}^{\varphi}}&=\mathds{1}(\phi=\varphi)\mathds{1}(j \in \mathbb{P}_{h}), \\
\frac{\partial \hat{\Lambda}_{i j}^{\phi}}{\partial q_{h}^{\varphi}}&=\mathbf{i}\mathds{1}(\phi=\varphi)\mathds{1}(j \in \mathbb{P}_{h}), \\
\frac{\partial \hat{v}_{j}^{\phi\phi}}{\partial p_{h}^{\varphi}}&= R_{jh}^{\phi\varphi}:=2 \sum_{(s, t) \in \mathbb{P}_{j \wedge h}} \operatorname{Re}\left(\bar{z}_{s t}^{\phi \varphi}\alpha^{\phi-\varphi} \right), \label{MLineargrad::p}\\
\frac{\partial \hat{v}_{j}^{\phi\phi}}{\partial q_{h}^{\varphi}} &= X_{jh}^{\phi\varphi}:=-2 \sum_{(s, t) \in \mathbb{P}_{j \wedge h}} \operatorname{Im}\left( \bar{z}_{s t}^{\phi \varphi}\alpha^{\phi-\varphi}\right), \label{MLineargrad::q}
\end{alignat}
\end{subequations}
where $\mathbb{P}_{j\wedge h}$ denotes the common part of the unique paths from bus $j$ and $h$ back to the root. Function $\mathds{1}(\phi=\varphi)$ is an indicator that equals $1$ if $\phi=\varphi$ and $0$ otherwise; $\mathds{1}(j \in \mathbb{P}_{h})=1$ if bus $j$ lies on the unique path from bus $h$ to the root, and $0$ otherwise. $\operatorname{Re}(\cdot)$ and $\operatorname{Im}(\cdot)$ denote the real and imaginary parts of a complex number, respectively.

% The above gradients have a similar structure to that in the single-phase case. The dynamics (\ref{MultiIter}) with simplified gradients (\ref{MLineargrad}) will converge to the solutions where the power network is exposed to the risk of voltage violation. This motivates us to propose a more accurate three-phase gradient evaluation method to prevent the voltage violation problem. 

Ignoring power loss (i.e., applying approximation Approx. 1) in the linearized BFM (\ref{linearmodel}) will induce an error in voltage prediction. To illustrate that error, we introduce an intermediate model that adds power loss back, while still making the balanced voltage approximation Approx. 2. This intermediate model is given by (\ref{disflow::v}) and (\ref{BTdisflow}) below:
\begin{subequations} \label{BTdisflow}
\begin{alignat}{2}
S_{i j}&=\gamma^{\Phi_{ij}} \operatorname{Diag}\left(\Lambda_{i j}\right),~  \forall (i,j)\in \mathcal{E}, \label{BTdisflow::v}\\
 \Lambda_{i j}&-\sum_{k: (j, k)\in \mathcal{E}} \Lambda_{j k}^{\Phi_j}=-s_{j}+\operatorname{diag}(z_{i j} \ell_{i j}), ~\forall j \in \mathcal{N}^+.\label{BTdisflow::PQ}
\end{alignat}
\end{subequations}
\begin{lemma} \label{lemma::Vsafe:three}
The difference between $\hat{v}_h$ predicted by (\ref{linearmodel}) and $v_h$ by (\ref{disflow::v}), (\ref{BTdisflow}) for all $h\in\mathcal{N}$ is:
\begin{eqnarray} \label{VestimateError::h}
\hat{v}_{h}-v_{h}=\sum_{(\zeta,\xi)\in \mathbb{P}_{h}}\left(M_{\zeta\xi}z_{\zeta\xi}^{H}+z_{\zeta\xi}M_{\zeta\xi}^{H}-z_{\zeta\xi} \ell_{\zeta\xi}z_{\zeta\xi}^{H}\right),
\end{eqnarray}
where $M_{\zeta\xi}$ depends on the power losses downstream of bus $\xi$ as:
\begin{eqnarray}
M_{\zeta\xi} & =&\gamma^{\Phi_{\zeta\xi}}\operatorname{Diag}(\operatorname{diag}(z_{\zeta\xi}\ell_{\zeta\xi}))  \nonumber\\
& &+\sum_{(\alpha,\beta)\in\operatorname{down}(\xi)}\gamma^{\Phi_{\alpha \beta}}\operatorname{Diag}(\operatorname{diag}(z_{\alpha\beta}\ell_{\alpha\beta})). \nonumber
\end{eqnarray}
\end{lemma} 
\begin{IEEEproof}
See Appendix \ref{appendix::1}.
\end{IEEEproof}

% \textit{Lemma \ref{lemma::Vsafe:three}} indicates the linearized model (\ref{linearmodel}) optimistically estimates the squared voltage magnitudes by $\sum_{(\zeta,\xi)\in \mathbbm{P}_{h}}\left(\left(M_{\zeta\xi}z_{\zeta\xi}^{H}+z_{\zeta\xi}M_{\zeta\xi}^{H}\right)-z_{\zeta\xi} \ell_{\zeta\xi}z_{\zeta\xi}^{H}\right)$. Similar voltage error analysis for the single-phase model was provided in \cite{liang2022hierarchical}. Therefore, under the power injections $(\boldsymbol{p}_{\Xi},\boldsymbol{q}_{\Xi})$ solved by (\ref{MultiIter}) with the linearized model (\ref{linearmodel}) and the simplified gradients (\ref{MLineargrad}), even though the model optimistically estimate the voltage to be safe, the actual voltage may already drop below their lower bounds. The need to prevent such voltage violation motivates us to develop an improved gradient evaluation method based on the accurate nonlinear three-phase unbalanced model, as elaborated below.

\textit{Lemma \ref{lemma::Vsafe:three}} reveals that the lossless model (\ref{linearmodel}) optimistically estimates the squared voltage magnitudes at all the buses, namely, even though the lossless model predicts the voltages to be safe, the actual voltages under the same power injections might already violate their safety limits. The need to prevent such voltage violation motivates us to develop an improved gradient evaluation method rather than simply using (\ref{MLineargrad}) from the lossless model (\ref{linearmodel}). A similar voltage error analysis was provided, and an improved gradient evaluation method was proposed in \cite{liang2022hierarchical} and its extended version on arXiv, for single-phase networks. This paper extends such analysis and method to three-phase unbalanced networks, as elaborated below.

\subsection{Improved Gradient Evaluation}

A major challenge for calculating the accurate gradients from BFM (\ref{distflow}) is the presence of the positive semidefinite (PSD) and rank-1 constraints (\ref{disflow::PSD})--(\ref{disflow::rank}). To deal with that, we relax (\ref{disflow::PSD})--(\ref{disflow::rank}) into equation (\ref{quadratic}) in \textit{Lemma \ref{lemma::2}} below. This relaxation extends the one in \cite{jha2021network}.
\begin{lemma} \label{lemma::2}
A necessary condition for PSD rank-1 constraints (\ref{disflow::PSD})--(\ref{disflow::rank}) is that for all $\phi,\psi, \eta \in \Phi_{ij}$, the following holds:
\begin{alignat}{2}\label{quadratic}
v_i^{\phi\phi} \ell_{ij}^{\psi\eta}=S_{ij}^{\phi\eta}\overline{S}_{ij}^{\phi\psi}.
\end{alignat}
\end{lemma}

The proof of \textit{Lemma \ref{lemma::2}} is straightforward from the linear dependence of any pair of rows (columns) in the rank-1 matrix.

As mentioned in Section \ref{subsec::motivation}, we replace (\ref{disflow::PQ}) with (\ref{BTdisflow}), where (\ref{BTdisflow::v}) utilizes (\ref{linearmodel::S}) based on the balanced voltage approximation Approx. 2, while (\ref{BTdisflow::PQ}) preserves the term $\operatorname{diag}(z_{ij}\ell_{ij})$ related to power loss. We call the modified model (\ref{disflow::v}), (\ref{BTdisflow}), (\ref{quadratic}) the BFM with balanced voltage approximation (BFM-BVA). Based on this model, the pertinent partial derivatives with respect to power injections $u_h^{\varphi}:=[p_h^{\varphi},q_h^{\varphi}]^{\top}$ at all buses $h\in\mathcal{N}$ and phases $\varphi \in \Phi_h$ are calculated as follows:
\begin{subequations}\label{Multigradient}
\begin{alignat}{2}
&\frac{\partial v_{j}^{\phi\phi}}{\partial u_h^{\varphi}}=\frac{\partial v_{i}^{\phi\phi}}{\partial u_h^{\varphi}}-\left(\sum_{\psi \in \Phi_{ij}}\frac{\partial S_{i j}^{\phi \psi}}{\partial u_{h}^{\varphi}} \overline{z}_{i j}^{\phi \psi}+\sum_{\psi \in \Phi_{ij}}z_{i j}^{\phi\psi} \frac{\partial \overline{S}_{i j}^{\phi \psi}}{\partial u_{h}^{\varphi}}\right) \nonumber\\
& +\sum_{\eta \in \Phi_{ij}}\sum_{\psi \in \Phi_{ij}}z_{i j}^{\phi \psi} \frac{\partial\ell_{i j}^{\psi\eta}}{\partial u_{h}^{\varphi}} z_{i j}^{\phi \eta}, ~\forall (i,j) \in \mathcal{E}, ~\forall \phi \in \Phi_{ij}, \label{Multigradient:v} \\
&\frac{\partial S_{i j}}{\partial u_{h}^{\varphi}}=\gamma^{\Phi_{i j}} \operatorname{Diag}\left(\frac{\partial \Lambda_{i j}}{\partial u_{h}^{\varphi}}\right),~  \forall (i,j)\in \mathcal{E}, \label{Multigradient:diag} \\
&\frac{\partial \Lambda_{ij}^{\phi}}{\partial u_h^{\varphi}}-\sum_{\psi \in \Phi_{ij}}z_{ij}^{\phi \psi}\frac{\partial \ell_{ij}^{\psi \phi}}{\partial u_h^{\varphi}}-\sum_{k: (j, k) \in \mathcal{E}}\frac{\partial \Lambda_{j k}^{\Phi_{j}\phi}}{\partial u_h^{\varphi}}=-\frac{\partial s_{j}^{\phi}}{\partial u_{h}^{\varphi}}, \nonumber\\
&\qquad \qquad \qquad \qquad \qquad ~\forall j \in \mathcal{N}^+,~\forall \phi \in \Phi_j,\\ 
&\frac{\partial\ell_{ij}^{\psi\eta}}{\partial u_{h}^{\varphi}}=\frac{\overline{S}_{ij}^{\phi \psi}}{v_{i}^{\phi\phi}}\frac{\partial S_{ij}^{\phi \eta}}{\partial u_{h}^{\varphi}}+\frac{S_{ij}^{\phi \eta}}{v_{i}^{\phi\phi}}\frac{\partial \overline{S}_{ij}^{\phi \psi}}{\partial u_{h}^{\varphi}}- \frac{\ell_{ij}^{\psi \eta}}{v_{i}^{\phi\phi}}\frac{\partial v_{i}^{\phi\phi}}{\partial u_{h}^{\varphi}}, \nonumber\\
&\qquad \qquad \qquad \qquad \qquad ~\forall (i,j)\in \mathcal{E}, ~\forall \psi,\eta,\phi\in \Phi_{ij}, \label{Multigradient::ell}
\end{alignat}
\end{subequations}
where $v_{i}^{\phi\phi}$, namely the $\phi$-th diagonal element in matrix $v_i$, is the squared voltage magnitude at bus $i$, phase $\phi$.

To decouple the complex interdependence of variables and partial derivatives in (\ref{Multigradient}), we simplify $\partial{\boldsymbol{\ell}}/\partial{\boldsymbol{u}}$ by replacing the partial derivatives therein with their approximations from the linearized model (\ref{linearmodel}):
\begin{alignat}{2} \label{Multiell}
\frac{\partial \hat{\ell}_{ij}^{\psi\eta}}{\partial u_{h}^{\varphi}}=\frac{\overline{S}_{ij}^{\phi \psi}}{v_{i}^{\phi\phi}}\frac{\partial \hat{S}_{ij}^{\phi \eta}}{\partial u_{h}^{\varphi}}+\frac{S_{ij}^{\phi \eta}}{v_{i}^{\phi\phi}}\frac{\partial \hat{\overline{S}}_{ij}^{\phi \psi}}{\partial u_{h}^{\varphi}}- \frac{\ell_{ij}^{\psi \eta}}{v_{i}^{\phi\phi}}\frac{\partial \hat{v}_{i}^{\phi\phi}}{\partial u_{h}^{\varphi}}.
\end{alignat}

By (\ref{MLineargrad}), we calculate $\partial{\hat{\boldsymbol{\ell}}}/\partial{\boldsymbol{u}}$ as:
\begin{subequations} \label{gradientell}
\begin{alignat}{2}
\frac{\partial\hat{\ell}_{ij}^{\psi\eta}}{\partial p_{h}^{\varphi}}=&-\frac{1}{v_{i}^{\phi\phi}}\left[\ell_{ij}^{\psi\eta}R_{ih}^{\phi\varphi}+ \left(\overline{S}_{ij}^{\phi \psi}\alpha^{\phi-\varphi}\mathbbm{1}(\eta=\varphi)\right.\right. \nonumber\\
&\left.\left.+S_{ij}^{\phi \eta}\overline{\alpha}^{(\phi-\varphi)}\mathbbm{1}(\psi=\varphi)\right)\mathbbm{1}(j \in \mathbbm{P}_{h})\right], \\
\frac{\partial\hat{\ell}_{ij}^{\psi\eta}}{\partial q_{h}^{\varphi}}=&-\frac{1}{v_{i}^{\phi\phi}}\left[\ell_{ij}^{\psi\eta}X_{ih}^{\phi\varphi}+ \left(\mathbf{i}\overline{S}_{ij}^{\phi \psi}\alpha^{\phi-\varphi}\mathbbm{1}(\eta=\varphi)\right.\right. \nonumber\\
&\left.\left.+\mathbf{i}S_{ij}^{\phi \eta}{\overline\alpha}^{(\phi-\varphi)}\mathbbm{1}(\psi=\varphi)\right)\mathbbm{1}(j \in \mathbbm{P}_{h})\right].
\end{alignat}
\end{subequations}

Replacing the partial derivatives on the RHS of (\ref{Multigradient:v}) with (\ref{MLineargrad}), (\ref{gradientell}) and (\ref{Multigradient:diag}), we get the following \emph{improved gradient evaluation}:
\begin{subequations} \label{MultiImpGrad}
\begin{alignat}{2}
&\frac{\partial \widetilde{v}_{j}^{\phi\phi}}{\partial p_h^{\varphi}}=\frac{\partial \hat{v}_{i}^{\phi\phi}}{\partial p_h^{\varphi}}-\left(\sum_{\psi \in \Phi_{ij}}\frac{\partial \hat{S}_{i j}^{\phi \psi}}{\partial p_{h}^{\varphi}} \overline{z}_{i j}^{\phi \psi}+\sum_{\psi \in \Phi_{ij}}z_{i j}^{\phi \psi} \frac{\partial \hat{\overline{S}}_{i j}^{\phi \psi}}{\partial p_{h}^{\varphi}}\right) \nonumber\\
&\quad\quad+\sum_{\psi \in \Phi_{ij}}\sum_{\eta \in \Phi_{ij}}z_{i j}^{\phi \psi} \frac{\partial \hat{\ell}_{i j}^{\psi\eta}}{\partial p_{h}^{\varphi}} z_{i j}^{\phi \eta} \nonumber \\
&=\left(1-\frac{1}{v_i^{\phi\phi}}\sum_{\psi \in \Phi_{ij}}\sum_{\eta\in\Phi_{ij}}\ell_{ij}^{\psi\eta}z_{ij}^{\phi\psi}\overline{z}_{ij}^{\phi\eta}\right)R_{ih}^{\phi\varphi}+\mathbbm{1}(j \in \mathbbm{P}_{h})\cdot \nonumber\\
&\left(2\operatorname{Re}(\overline{z}_{ij}^{\phi\varphi}\alpha^{\phi-\varphi})-\frac{2}{v_{i}^{\phi\phi}}\sum_{\psi\in \Phi_{ij}}\operatorname{Re}(\alpha^{\phi-\varphi}\overline{S}_{ij}^{\phi \psi}z_{ij}^{\phi\psi}\overline{z}_{ij}^{\phi\varphi})\right), \label{MultiImpGrad::p}\\
&\frac{\partial \widetilde{v}_{j}^{\phi\phi}}{\partial q_h^{\varphi}}=\frac{\partial \hat{v}_{i}^{\phi\phi}}{\partial q_h^{\varphi}}-\left(\sum_{\psi \in \Phi_{ij}}\frac{\partial \hat{S}_{i j}^{\phi \psi}}{\partial q_{h}^{\varphi}} \overline{z}_{i j}^{\phi \psi}+\sum_{\psi \in \Phi_{ij}}z_{i j}^{\phi \psi} \frac{\partial \hat{\overline{S}}_{i j}^{\phi \psi}}{\partial q_{h}^{\varphi}}\right) \nonumber\\
&\quad\quad+\sum_{\psi \in \Phi_{ij}}\sum_{\eta \in \Phi_{ij}}z_{i j}^{\phi \psi} \frac{\partial \hat{\ell}_{i j}^{\psi\eta}}{\partial q_{h}^{\varphi}} z_{i j}^{\phi \eta} \nonumber \\
&=\left(1-\frac{1}{v_i^{\phi\phi}}\sum_{\psi \in \Phi_{ij}}\sum_{\eta\in\Phi_{ij}}\ell_{ij}^{\psi\eta}z_{ij}^{\phi\psi}\overline{z}_{ij}^{\phi\eta}\right)X_{ih}^{\phi\varphi}+\mathbbm{1}(j \in \mathbbm{P}_{h})\cdot \nonumber\\
&\left(-2\operatorname{Im}(\overline{z}_{ij}^{\phi\varphi}\alpha^{\phi-\varphi})+\frac{2}{v_{i}^{\phi\phi}}\sum_{\psi\in \Phi_{ij}}\operatorname{Im}(\alpha^{\phi-\varphi}\overline{S}_{ij}^{\phi \psi}z_{ij}^{\phi\psi}\overline{z}_{ij}^{\phi\varphi})\right). \label{MultiImpGrad::q}
\end{alignat}
\end{subequations}

The major difference of the improved gradient (\ref{MultiImpGrad::p}) from the previous approximate gradient (\ref{MLineargrad::p}) is the following (similar to that between (\ref{MultiImpGrad::q}) and (\ref{MLineargrad::q})). By (\ref{MLineargrad::p}), we have $R_{jh}^{\phi\varphi}=R_{ih}^{\phi\varphi}+2\operatorname{Re}(\overline{z}_{ij}^{\phi\varphi}\alpha^{\phi-\varphi})\cdot\mathds{1}(j \in \mathbb{P}_{h})$, which can convert (\ref{MultiImpGrad::p}) into:
\begin{alignat}{2} \label{MultiImpGrad:rewrite}
&\frac{\partial \widetilde{v}_{j}^{\phi\phi}}{\partial p_h^{\varphi}}=R_{jh}^{\phi\varphi}-\frac{R_{ih}^{\phi\varphi}}{v_i^{\phi\phi}}\sum_{\psi \in \Phi_{ij}}\sum_{\eta\in\Phi_{ij}}\ell_{ij}^{\psi\eta}z_{ij}^{\phi\psi}\overline{z}_{ij}^{\phi\eta} \nonumber \\
&-\left(\frac{2}{v_{i}^{\phi\phi}}\sum_{\psi\in \Phi_{ij}}\operatorname{Re}(\alpha^{\phi-\varphi}\overline{S}_{ij}^{\phi \psi}z_{ij}^{\phi\psi}\overline{z}_{ij}^{\phi\varphi})\right)\cdot \mathbbm{1}(j \in \mathbbm{P}_{h}).
\end{alignat}
The first term on the RHS of (\ref{MultiImpGrad:rewrite}) is exactly (\ref{MLineargrad::p}), while the second and the third terms partly compensate for the error due to ignoring power loss. Compared to the improved gradient evaluation in \cite{liang2022hierarchical} for single-phase networks, the proposed evaluation (\ref{MultiImpGrad}) considers mutual impedance between different phases in three-phase networks. 

\section{Improved Hierarchical OPF Algorithm}
\label{sec:multi:algorithm}

Based on the improved gradients (\ref{MultiImpGrad}) and the primal-dual algorithm (\ref{MultiIter}), we design a scalable algorithm to solve the OPF problem (\ref{MultiOPF}). This algorithm shares a hierarchical structure similar to that in \cite{zhou2019accelerated}, but is more reliable by using the proposed improved gradient evaluation, which can alleviate the voltage violation induced by ignoring power loss.

The tree network $\mathcal{T}:=\{\mathcal{N}^{+},\mathcal{E}\}$ is divided into subtrees $\mathcal{T}_{k}=\{\mathcal{N}_{k},\mathcal{E}_{k}\}$ indexed by $k\in \mathcal{K}=\{1,\dots,K\}$ and a set of buses $\mathcal{N}_{0}$ that are not clustered into any subtree. Let $n_{k}^{0}$ denote the root bus of subtree $k$, which is the bus in $\mathcal{N}_{k}$ that is nearest to the root bus of the whole network. The division of subtrees and unclustered buses may not be unique, but we assume it always satisfies the following conditions:
\begin{assumption}\label{ass:nooverlap}
All the subtrees are non-overlapping, i.e., $\mathcal{N}_{k_{1}}\cap \mathcal{N}_{k_2}=\emptyset$ for any $k_1,k_2 \in \mathcal{K},k_1 \neq k_2
$.
\end{assumption}
\begin{assumption} \label{ass:notonpath}
For any subtree root bus $n_{k}^{0}, k \in \mathcal{K}$, or any unclustered bus $n \in \mathcal{N}_{0}$, its path to the network root bus only goes through a subset of buses in $\mathcal{N}_{0}$, but not any bus in another subtree. 
\end{assumption}

\subsection{Three-Phase Hierarchical Algorithm}
\begin{figure} 
    \centering
    \includegraphics[width = 0.49\columnwidth]{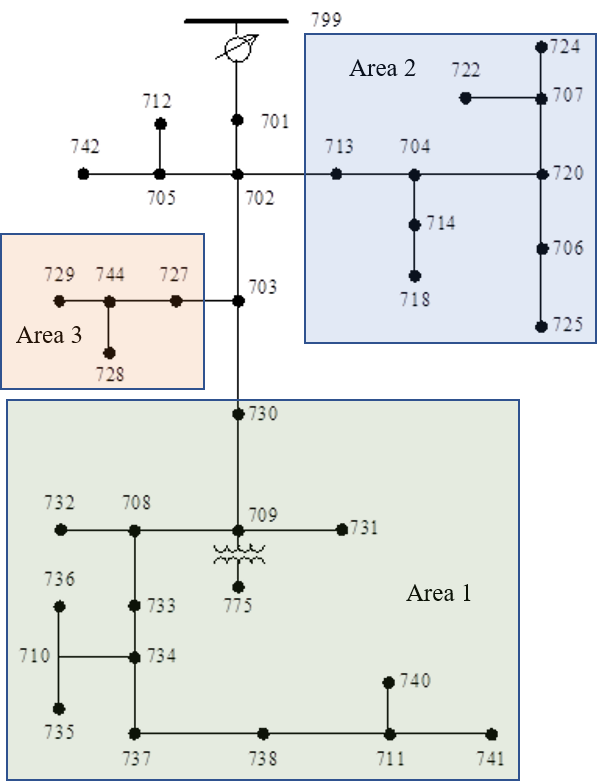}
    \hfil
    \includegraphics[width = 0.49\columnwidth]{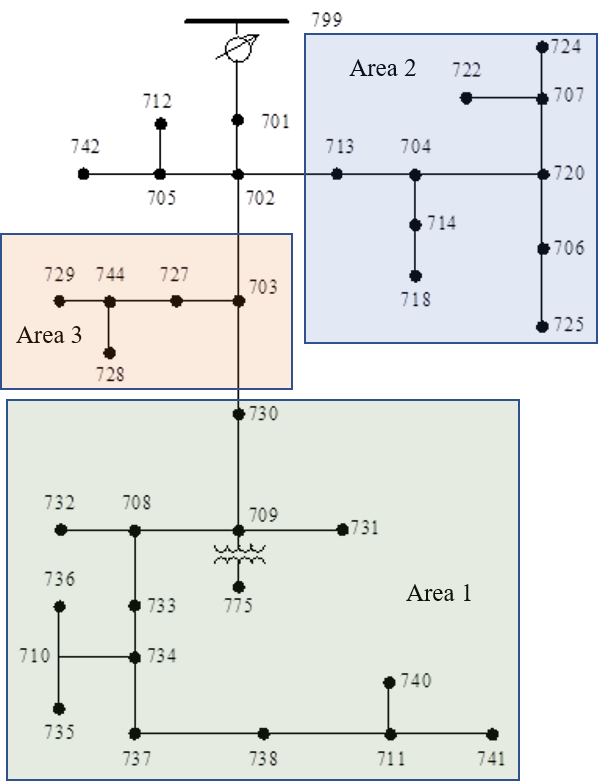}
    \caption{The clustering of the IEEE $37$-bus network on the left satisfies Assumption \ref{ass:notonpath}. The one on the right does not, because the path from a subtree root bus 730 to the network root 799 goes through another subtree.}
    \label{ieee37}
\end{figure}

The left subfigure of Figure \ref{ieee37} shows a clustering of the IEEE $37$-bus network that satisfies both Assumptions \ref{ass:nooverlap} and \ref{ass:notonpath}, while the one in the right subfigure satisfies Assumption \ref{ass:nooverlap} only, but not Assumption \ref{ass:notonpath}. How to optimally cluster a tree network under these two conditions is beyond the scope of this paper but is one of the interests for our future research.

The above settings facilitate a scalable hierarchical management of the power network. Each subtree $k\in \mathcal{K}$ is represented by its root bus $n_k^0$ and managed by a regional controller (RC) that is also indexed by $k$. The subtree root buses together with the unclustered buses form a reduced tree, which is managed by a central coordinator (CC). Both the CC and RCs maintain two-way communication with their respective buses being managed. The hierarchical structure of CC-RC-subtree was proposed in \cite{zhou2019accelerated} and now we adapt it to incorporate the improved gradient evaluation (\ref{MultiImpGrad}). In particular, the primal-dual algorithm (\ref{MultiIter}) can be implemented locally and separately at each bus, except for the term $\frac{\partial \boldsymbol{v}}{\partial \boldsymbol{u}}^{\top}(\overline{\boldsymbol{\mu}}-\underline{\boldsymbol{\mu}})$ that in principle couples the entire network across buses. We now explain how to leverage the hierarchical structure to accelerate computation of this coupling term. We only elaborate the terms associated with active power injection $p_{h}^{\varphi}$ at phase $\varphi$ of bus $h \in \mathcal{N}$, as those related to reactive power $q_{h}^{\varphi}$ are calculated similarly. There are two cases to consider:
% \begin{figure*}
%   \centering
%     \includegraphics[width=1.55\columnwidth]{figures/CC_RCs.png}
%     \caption{The clustering of IEEE 123-node test feeder in our experiments.}
%     \label{fig:ieee123}
% \end{figure*}

\textit{Case 1:} If $h \in \mathcal{N}_{k}$ is in subtree $k \in \mathcal{K}$, we have:
\begin{alignat}{2}
 &\sum_{j \in \mathcal{N}}\sum_{\phi \in \Phi_{j}} \frac{\partial \widetilde{v}_{j}^{\phi\phi}}{\partial p_{h}^{\varphi}}\left(\overline{\mu}_{j}^{\phi}-\underline{\mu}_{j}^{\phi}\right)=\sum_{\phi \in \Phi_{0}}\sum_{j \in \mathcal{N}_{k}^{\phi}}\frac{\partial \widetilde{v}_{j}^{\phi\phi}}{\partial p_{h}^{\varphi}}\left(\overline{\mu}_{j}^{\phi}-\underline{\mu}_{j}^{\phi}\right) \nonumber\\
 &\qquad +\sum_{\phi\in \Phi_{0}}\bigg[\sum_{k' \in \mathcal{K}\backslash\{k\}}~\sum_{j \in \mathcal{N}_{k'}^{\phi}}\frac{\partial \widetilde{v}_{j}^{\phi\phi}}{\partial p_{h}^{\varphi}}\left(\overline{\mu}_{j}^{\phi}-\underline{\mu}_{j}^{\phi}\right) 
 \nonumber \\ 
 & \qquad +\sum_{j \in \mathcal{N}_{0}^{\phi}}\frac{\partial \widetilde{v}_{j}^{\phi\phi}}{\partial p_{h}^{\varphi}}\left(\overline{\mu}_{j}^{\phi}-\underline{\mu}_{j}^{\phi}\right)\bigg],\label{eq:decomposed_dvdp}
 \end{alignat}
where $\Phi_{0}$ denotes the phase set of the whole network, and $\mathcal{N}_{k}^{\phi}$ (or $\mathcal{N}_{0}^{\phi}$) denotes the subset of buses in subtree $k$ (or the subset of unclustered buses) that have phase $\phi$.

The first term on the RHS of (\ref{eq:decomposed_dvdp}) sums over the buses in the same subtree $k$ as bus $h$, which can be calculated by RC $k$ using (\ref{MultiImpGrad::p}). Specifically, from a power-flow solution in each iteration, RC $k$ gets the current, power flow, and voltage $(\ell_{ij}, S_{ij}, v_{i})$ associated with each bus $j\in \mathcal{N}_{k}^{\phi}$, while it also knows impedance $z_{ij}$ and whether $j \in \mathbb{P}_{h}$. Upon gathering such information within subtree $k$, RC $k$ conducts the calculation (\ref{MultiImpGrad::p}), and then sends the result to bus $h$.

The third term on the RHS of (\ref{eq:decomposed_dvdp}) sums over unclustered buses $j \in \mathcal{N}_{0}^{\phi}$. Since all such buses are located in the reduced tree, this term can be calculated by the CC using (\ref{MultiImpGrad::p}). Moreover, from the outside of subtree $k$, all of its buses can be represented by its root bus $n_{k}^{0}$. In particular, $R_{ih}^{\phi\varphi}=R_{in_{k}^{0}}^{\phi\varphi}$ and $\mathds{1}(j\in\mathbb{P}_{h}) = \mathds{1}(j\in\mathbb{P}_{n_k^0})$ are known to the CC for the calculation of (\ref{MultiImpGrad::p}). 

The second term on the RHS of (\ref{eq:decomposed_dvdp}) sums over all buses $j$ in subtrees other than subtree $k$ that hosts bus $h$. By Assumption \ref{ass:notonpath}, there must be $j \notin \mathbb{P}_{h}$. Moreover, there is $R_{ih}^{\phi\varphi}=R_{n_k^{0} n_{k'}^{0}}^{\phi\varphi}$ for all $h \in \mathcal{N}_{k},j \in \mathcal{N}_{k'}^{\phi}$, and $i$ being the unique upstream bus of $j$. This, by (\ref{MultiImpGrad::p}), simplifies this term into:
\begin{alignat}{2} \label{NkNd}
&\sum_{k' \in \mathcal{K}\backslash\{k\}}~\sum_{j \in \mathcal{N}_{k'}^{\phi}}\frac{\partial \widetilde{v}_{j}^{\phi\phi}}{\partial p_{h}^{\varphi}}\left(\overline{\mu}_{j}^{\phi}-\underline{\mu}_{j}^{\phi}\right)=\sum_{k' \in \mathcal{K}\backslash\{k\}} R_{n_k^0 n_{k'}^0}^{\phi\varphi}\cdot \nonumber\\
&\sum_{j \in \mathcal{N}_{k'}^{\phi}}\left(1-\frac{1}{v_i^{\phi\phi}}\sum_{\psi \in \Phi_{ij}}\sum_{\eta\in\Phi_{ij}}\ell_{ij}^{\psi\eta}z_{ij}^{\phi\psi}\overline{z}_{ij}^{\phi\eta}\right)\left(\overline{\mu}_{j}^{\phi}-\underline{\mu}_{j}^{\phi}\right).
\end{alignat}
In (\ref{NkNd}), the summation over $j \in \mathcal{N}_{k'}^{\phi}$ is calculated by RC $k'$ and then uploaded to the CC. The CC receives such results from all RCs $k' \in \mathcal{K}\backslash\{k\}$ and adds them up after weighting by $R_{n_k^{0} n_{k'}^{0}}^{\phi\varphi}$. The CC finally adds the second and third terms on the RHS of (\ref{eq:decomposed_dvdp}) over all phases $\phi \in \Phi_0$, and then sends the result to each RC $k$. Then RC $k$ broadcasts that result to all the buses in subtree $k$, which then calculates (\ref{eq:decomposed_dvdp}) and (\ref{MultiIter:p}).

\textit{Case 2:} If $h \in \mathcal{N}_{0}$ is an unclustered bus, we have:
\begin{alignat}{2}
&\sum_{j \in \mathcal{N}}\sum_{\phi \in \Phi_{j}}\frac{\partial \widetilde{v}_{j}^{\phi\phi}}{\partial p_{h}^{\varphi}}\left(\overline{\mu}_{j}^{\phi}-\underline{\mu}_{j}^{\phi}\right)=\sum_{\phi \in \Phi_{0}}\bigg[\sum_{j \in \mathcal{N}_{0}^{\phi}}\frac{\partial \widetilde{v}_{j}^{\phi\phi}}{\partial p_{h}^{\varphi}}\left(\overline{\mu}_{j}^{\phi}-\underline{\mu}_{j}^{\phi}\right)\nonumber\\
&\qquad +\sum_{k \in \mathcal{K}}~\sum_{j \in \mathcal{N}_{k}^{\phi}}\frac{\partial \widetilde{v}_{j}^{\phi\phi}}{\partial p_{h}^{\varphi}}\left(\overline{\mu}_{j}^{\phi}-\underline{\mu}_{j}^{\phi}\right)\bigg]. \label{eq:decomposed_dvdp_unclustered}
\end{alignat}

\begin{figure}
    \centering
    \includegraphics[width=0.90\columnwidth]{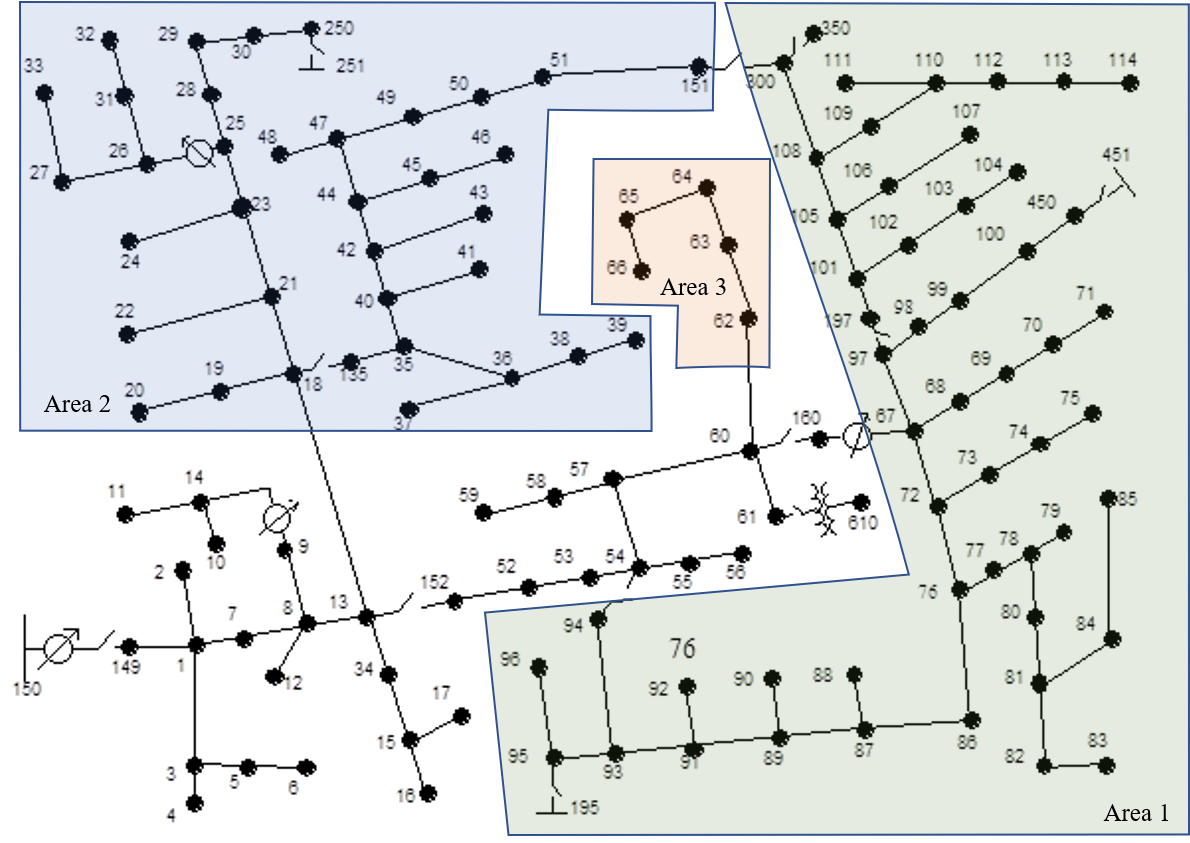}
    \caption{The clustering of IEEE $123$-bus test feeder in our experiments.}
    \label{fig:ieee123}
\end{figure}

The idea to calculate (\ref{eq:decomposed_dvdp_unclustered}) through information exchange between bus $h$, all RCs $k\in\mathcal{K}$, and the CC is the same as that for the second and third terms in (\ref{eq:decomposed_dvdp}). The only difference is that there is no longer a subtree managed by bus $h$, so the process becomes even simpler.

To summarize, the computation of the key term $\frac{\partial \boldsymbol{v}}{\partial \boldsymbol{u}}^{\top}(\overline{\boldsymbol{\mu}}-\underline{\boldsymbol{\mu}})$ in the primal-dual algorithm (\ref{MultiIter}) can be performed through the coordination of the CC, RCs, and unclustered buses. They form a hierarchical structure, which inspires us to design Algorithm \ref{algCCRC} to accelerate such computation. Compared to the conventional centralized primal-dual gradient method, Algorithm \ref{algCCRC} can reduce computational complexity as can be analyzed in the same way as \cite{zhou2019accelerated}. However, compared to \cite{zhou2019accelerated} based on a lossless model, our improved gradient evaluation considers the effect of power loss on voltage estimation and thus more reliably guarantees voltage safety.

\begin{algorithm*}
\caption{Hierarchical OPF Algorithm For Three-Phase Unbalanced Distribution Networks} \label{algCCRC}
\begin{algorithmic}[1]
% \State \algorithmicgiven ~Initial feasible point $\bar{\alpha}_0$, $\bar{\beta}_0$, stopping criterion $\delta > 0$, $k := 0$.
% %
\Repeat
\State At time step $t$, every bus $h\in \mathcal{N}$ updates the power injection $(p_h^{\varphi}(t+1),q_h^{\varphi}(t+1))$ into the power network for all phases $\varphi \in \Phi_{h}$ by:
\begin{subequations}
\begin{alignat}{2}
p_{h}^{\varphi}(t+1)=\left[p_{h}^{\varphi}(t)-\sigma_u\left(\frac{\partial f_{h}^{\varphi}\left({p}_h^{\varphi}(t), {q}_h^{\varphi}(t)\right)}{\partial p_{h}^{\varphi}}+\alpha_{h}^{\varphi}(t)\right)  \right]_{\mathcal{Y}_{h}^{\varphi}} \nonumber\\
q_{h}^{\varphi}(t+1)=\left[q_{h}^{\varphi}(t)-\sigma_u\left(\frac{\partial f_{h}^{\varphi}\left({p}_h^{\varphi}(t), {q}_h^{\varphi}(t)\right)}{\partial q_{h}^{\varphi}}+\beta_{h}^{\varphi}(t)\right)  \right]_{\mathcal{Y}_{h}^{\varphi}} \nonumber
\end{alignat}
\end{subequations}
where $\alpha_{h}^{\varphi}(t)=\quad \sum_{j \in \mathcal{N}}\sum_{\phi\in \Phi_{j}}\frac{\partial \widetilde{v}_{j}^{\phi\phi}\left(\boldsymbol{u}(t)\right)}{\partial p_{h}^{\varphi}}\cdot\left(\overline{\mu}_{j}^{\varphi}(t)-\underline{\mu}_{j}^{\varphi}(t)\right)$ is given by \eqref{eq:decomposed_dvdp} if $h\in\mathcal{N}_k$ is in subtree $k$, or \eqref{eq:decomposed_dvdp_unclustered} if $h\in\mathcal{N}_0$ is unclustered; and $\beta_{h}^{\varphi}(t)=\quad \sum_{j \in \mathcal{N}}\sum_{\phi \in \Phi_{j}}\frac{\partial \widetilde{v}_{j}^{\phi\phi}\left(\boldsymbol{u}(t)\right)}{\partial q_{h}^{\varphi}}\cdot\left(\overline{\mu}_{j}^{\varphi}(t)-\underline{\mu}_{j}^{\varphi}(t)\right)$. 

\State After a power flow is solved using the updated power injections, every bus $h\in\mathcal{N}$ uses its local voltage $v_h^{\varphi\varphi}$ for all phases $\varphi \in \Phi_h$ to update its dual variables:
\begin{alignat}{2}
\underline{\mu}_{h}^{\varphi}(t+1)=\left[\underline{\mu}_{h}^{\varphi}(t)+\sigma_{\mu}\big(\underline{v}_{h}^{\varphi}-v_{h}^{\varphi\varphi}(t)-\epsilon \underline{\mu}_{h}^{\varphi}(t)\big)\right]_{+}, \quad \overline{\mu}_{h}^{\varphi}(t+1)=\Big[\overline{\mu}_{h}^{\varphi}(t)+\sigma_{\mu}\big( v_{h}^{\varphi\varphi}(t)-\overline{v}_{h}^{\varphi}-\epsilon \overline{\mu}_{h}^{\varphi}(t)\big)\Big]_{+}. \nonumber
\end{alignat}
\State In every RC $k\in\mathcal{K}$, every bus $j\in \mathcal{N}_{k}$ sends local information $(v_{i},\ell_{ij},z_{ij})$ to RC $k$. Then, every RC $k\in \mathcal{K}$ calculates the following weighted sum of dual variables for all phases $\phi \in \Phi_{0}$, and sends it to the CC:
\begin{alignat}{2}
\sum_{j \in \mathcal{N}_{k}^{\phi}}\left(1-\frac{1}{v_i^{\phi\phi}}\sum_{\psi \in \Phi_{ij}}\sum_{\eta\in\Phi_{ij}}\ell_{ij}^{\psi\eta}z_{ij}^{\phi\psi}\overline{z}_{ij}^{\phi\eta}\right)\cdot\left(\overline{\mu}_{j}^{\phi}(t+1)-\underline{\mu}_{j}^{\phi}(t+1)\right), ~\forall \phi \in \Phi_{0}. \nonumber
\end{alignat}

\State The CC computes the second term and third term on the RHS of (\ref{eq:decomposed_dvdp}), and then adds them over all phases $\phi \in \Phi_0$ for each destination bus $h \in \mathcal{N}_k$ in subtree $k$, and finally sends the result to corresponding RC $k$ for all $k \in \mathcal{K}$. The CC also computes the first term and second term on the RHS of (\ref{eq:decomposed_dvdp_unclustered}), and then adds them over all phases $\phi \in \Phi_0$ for each unclustered destination bus $h\in \mathcal{N}_{0}$, and finally sends the result to each bus $h\in \mathcal{N}_{0}$. The terms needed for computing $\beta_h^{\varphi}(t+1)$ are calculated and sent by the CC similarly.
% \begin{alignat}{2}
% &s_{k'}^{\varphi out}=\sum_{\phi\in \Phi}\bigg(\sum_{k'\in\mathcal{K}\backslash\{k\}}R_{n_k^0 n_{k'}^0}^{\phi\varphi} \cdot B^{\phi}_{k}+U_{k}\bigg), \qquad s_{h}^{\varphi}=\sum_{\phi\in \Phi}\bigg(\sum_{k'\in\mathcal{K}} R_{h n_{k'}^0}^{\phi\varphi} \cdot B^{\phi}_{k'} +U_0 \bigg).
% \end{alignat}

\State For the bus $h\in \mathcal{N}_k$ in subtree $k$, RC $k$ computes the first term on the RHS of (\ref{eq:decomposed_dvdp}), adds it with the result received from the CC in step 5, and sends the result $\alpha_h^{\varphi}(t+1)$ to each destination bus $h$. For the unclustered bus $h \in \mathcal{N}_0$, it already obtains $\alpha_h^{\varphi}(t+1)$ in step 5 as (\ref{eq:decomposed_dvdp_unclustered}). The term $\beta_h^{\varphi}(t+1)$ is obtained similarly.

\Until $\left\|\boldsymbol{u}(t+1)-\boldsymbol{u}(t)\right\|_2 < \lambda$ for some preset threshold $\lambda>0$, or a maximum number of iterations is reached. 

\State Otherwise $t\leftarrow (t+1)$ and go back to step 2. 
\end{algorithmic}
\end{algorithm*}

\subsection{Convergence Analysis}

In this part, we further justify the significance of improved gradient evaluation by analyzing the convergence of Algorithm \ref{algCCRC} to the local optimal point(s) of the actual nonconvex problem (\ref{MultiOPF}) with improved gradient evaluation.

To facilitate our analysis, we rewrite Lagrangian (\ref{MultiRLagarangian}) in a compact form:
\begin{alignat}{2} \label{MultiLagran}
\mathcal{L}_{\epsilon}(\boldsymbol{u},\boldsymbol{\mu})=f(\boldsymbol{u})+\boldsymbol{\mu}^{\top}(A\boldsymbol{v}(\boldsymbol{u})-d)-\frac{\epsilon}{2}\|\boldsymbol{\mu}\|^2_2,
\end{alignat}
where $A$ and $d$ are the transformation matrices to make (\ref{MultiLagran}) aligned with (\ref{MultiRLagarangian}), i.e., $A=[\operatorname{diag}(\boldsymbol{1}_N), \operatorname{diag}(-\boldsymbol{1}_{N})]^{\top}$, and $d=[\underline{\boldsymbol{v}}^{\top},-\overline{\boldsymbol{v}}^{\top}]^{\top}$.

The dynamics (\ref{MultiIter}) with improved gradient evaluation (\ref{MultiImpGrad}) to approach a saddle point of (\ref{MultiLagran}) take the form:
\begin{subequations} \label{MultiIter:rewrite}
\begin{alignat}{2}
&\boldsymbol{u}(t+1)=\bigg[\boldsymbol{u}(t)-\sigma\bigg(\nabla f(\boldsymbol{u}(t))+\frac{\partial \widetilde{\boldsymbol{v}}(\boldsymbol{u}(t))}{\partial \boldsymbol{u}}^{\top}A^{\top}\boldsymbol{\mu}(t)\bigg)\bigg]_{\mathcal{Y}}, \\
&\boldsymbol{\mu}(t+1)=\bigg[\boldsymbol{\mu}(t)+\nu\sigma\left(A\boldsymbol{v}(\boldsymbol{u}(t))-d-\epsilon \boldsymbol{\mu}(t)\right)\bigg]_+.
\end{alignat}
\end{subequations}

A primal-dual local optimal point $(\boldsymbol{u}^*,\boldsymbol{\mu}^*)$ of the nonconvex OPF (\ref{MultiOPF}) satisfies the KKT condition:
\begin{subequations}\label{KKT}
\begin{alignat}{2}
\boldsymbol{u}^*\in \mathcal{Y}, \quad \boldsymbol{\mu}^*\geq 0,\\
\nabla f(\boldsymbol{u}^*)+\frac{\partial \boldsymbol{v}(\boldsymbol{u}^*)}{\partial \boldsymbol{u}}^{\top}A^{\top}\boldsymbol{\mu}^* \in -\mathbbm{N}_{\mathcal{Y}}(\boldsymbol{u^*}), \\
A\boldsymbol{v}(\boldsymbol{u}^*)-d\geq 0, \quad \boldsymbol{\mu}^{*\top}(A\boldsymbol{v}(\boldsymbol{u}^*)-d)=0,
\end{alignat}
\end{subequations}
where $\mathbbm{N}_{\mathcal{Y}}(\boldsymbol{u^*})$ is the normal cone of convex set $\mathcal{Y}$ at $\boldsymbol{u^*}\in\mathcal{Y}$, defined as $\{\boldsymbol{y}\in\mathbb{R}^{2(N-1)}:\boldsymbol{y}^{\top}(\boldsymbol{u}-\boldsymbol{u}^*)\leq 0,~\forall \boldsymbol{u} \in \mathcal{Y}\}$.

The conditions, under which the KKT point above is globally optimal, were discussed in \cite{farivar2013branch,gan2016online}. However, their analyses were based on single-phase networks. The conditions to guarantee global optimality of (\ref{KKT}) in three-phase networks remain an open question for our future research.

\begin{assumption}\label{ass:boundgradient}
The discrepancy between gradient evaluation (\ref{MultiImpGrad}) and the actual gradient given by the nonlinear model (\ref{distflow}) is bounded, i.e., there exists positive constant $e_1$ such that for all power injections $\boldsymbol{u}:=(\boldsymbol{p},\boldsymbol{q})\in\mathcal{Y}$:
\begin{alignat}{2}\label{def:boundgradient}
\left\|\frac{\partial \widetilde{\boldsymbol{v}}(\boldsymbol{u}(t))}{\partial \boldsymbol{u}}-\frac{\partial \boldsymbol{v}(\boldsymbol{u}(t))}{\partial \boldsymbol{u}}\right\|\leq e_1.
\end{alignat}
\end{assumption}

Assumption \ref{ass:boundgradient} generally holds by our inference in Section \ref{sec::improved}. We define:
\begin{subequations} \label{def::nonlinear}
\begin{alignat}{2}
& L_{\boldsymbol{v}}(\delta)=\sup_{\tau:\|\tau\|\leq \delta}\left\|\frac{\partial \boldsymbol{v}(\boldsymbol{u}^*+\tau)}{\partial \boldsymbol{u}}\right\|, \quad M_{\boldsymbol{\mu}}= \left\|\boldsymbol{\mu}^*\right\|, \\
&M_{\boldsymbol{v}}(\delta)=\sup_{\tau:\|\tau\|\leq \delta}\frac{\left\|\boldsymbol{v}(\boldsymbol{u}^*+\tau)-\boldsymbol{v}(\boldsymbol{u}^*)-\frac{\partial \boldsymbol{v}(\boldsymbol{u}^*+\tau)}{\partial \boldsymbol{u}}^{\top}\tau\right\|}{\|\tau\|^2}.
\end{alignat}
\end{subequations}

We provide a general explanation to the above definitions. $L_{\boldsymbol{v}}(\delta)$ represents the Lipschitz constant and $M_{\boldsymbol{v}}(\delta)$ characterizes the nonlinearity of function $\boldsymbol{v}(\boldsymbol{u})$ within a neighborhood of radius $\delta$ around the primal optimizer $\boldsymbol{u}^*$. See \cite{tang2018feedback,tang2022running} and references therein. Our difference lies in analyzing the impact of approximate gradient on the convergence of dynamics (\ref{MultiIter:rewrite}), while references \cite{tang2018feedback,tang2022running} assumed an accurate gradient.

\begin{figure*}
\begin{minipage}[h]{0.49\linewidth}
\centering
\includegraphics[width=1.00\columnwidth]{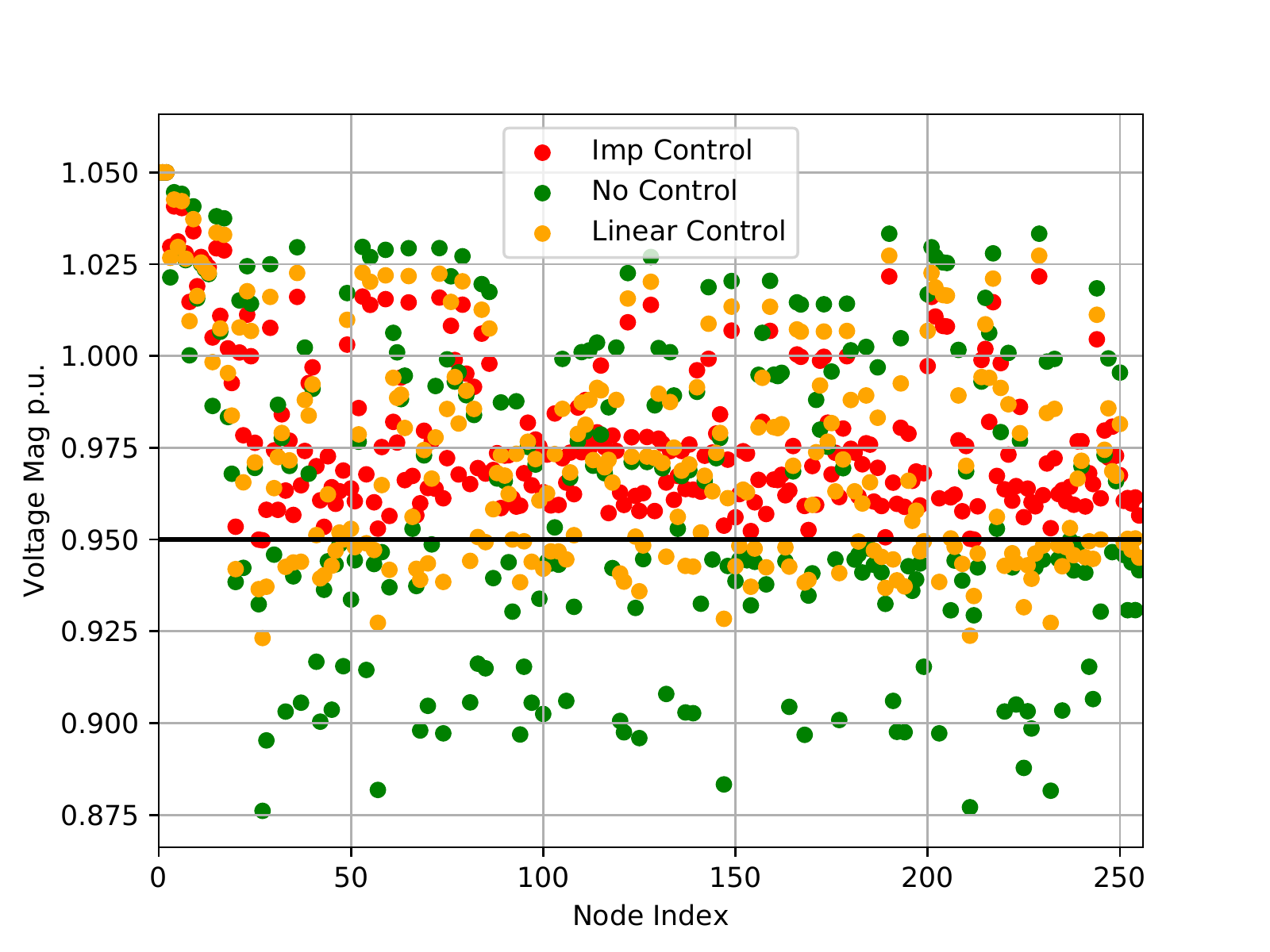}
% \caption{The voltage change at nodes of three different subtrees in the 4518-node network.}
% \label{fig:voltage_convergence_11000}
\end{minipage}
\begin{minipage}[h]{0.49\linewidth}
\centering
\includegraphics[width=1.00\columnwidth]{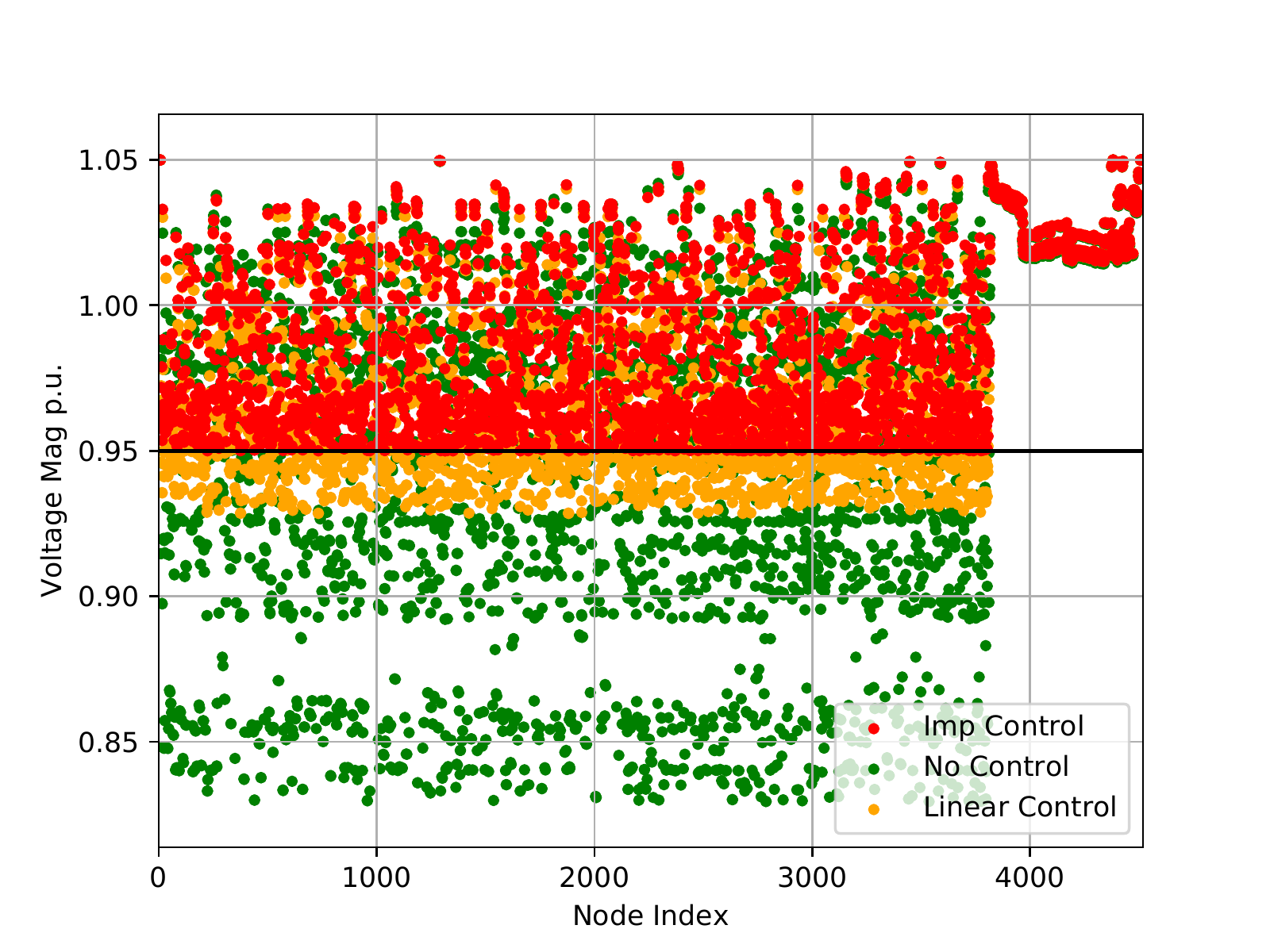}
\end{minipage}
\caption{The nodal voltages in the IEEE $123$-bus feeder (left) and the $4,518$-node feeder (right), in three cases: ``no control'' (using the nominal power injections), ``linear control \cite{zhou2019accelerated}'' (solving OPF based on the lossless linear model), and ``improved control'' (using Algorithm 1 with improved gradient evaluation).}
\label{fig:voltage_scatter}
\end{figure*}

We define: 
\begin{alignat}{2} \label{def:Delta}
\Delta=\inf_{\tau:\|\tau\|\leq\delta}\int_0^1\nabla_{\boldsymbol{u}\boldsymbol{u}}^2\mathcal{L}_{\epsilon}\left(\boldsymbol{u}^*+\theta\tau,\boldsymbol{\mu}^*\right)d\theta,
\end{alignat}
and
\begin{alignat}{2} \label{convergeSpeed}
&\rho(\delta,\sigma,e_1)= \bigg[\max\{(1-\sigma\Delta)^2,(1-\nu\sigma\epsilon)^2\}+\sigma\|A\|\nu^{\frac{1}{2}} \nonumber\\
&\cdot\left[(1-\nu\sigma\epsilon)(M_{\boldsymbol{v}}(\delta)\delta+e_1)+\sigma|\Delta-\nu\epsilon|(L_{\boldsymbol{v}}(\delta)+e_1)\right] \nonumber \\
&\quad +\nu\sigma^2\|A\|^2(L_{\boldsymbol{v}}(\delta)+e_1)^2\bigg]^\frac{1}{2}.
\end{alignat}

To measure the convergence, we define the norm:
\begin{alignat}{2} \nonumber
\left\|\boldsymbol{z}\right\|_{\nu}^2:=\left\|\boldsymbol{u}\right\|^2+\nu^{-1}\left\|\boldsymbol{\mu}\right\|^2.
\end{alignat}

We then have the following theorem regarding convergence of dynamic (\ref{MultiIter:rewrite}).
\begin{theorem}\label{theom:1}
Suppose there exists $\sigma>0$, $\nu>0$, $\epsilon>0$, and $\delta>0$ such that:
$$\rho(\delta,\sigma,e_1)<1.$$
And the initialization satisfies:
\begin{alignat}{2}
\left\|\boldsymbol{z}(1)-\boldsymbol{z}^*\right\|_{\nu}\leq \delta. \nonumber
\end{alignat}

The dynamics (\ref{MultiIter:rewrite}) converge to a neighborhood of $\boldsymbol{z}^*$ as:
\begin{alignat}{2} \label{convergeRegion}
\limsup_{t \rightarrow \infty}\left\|\boldsymbol{z}(t)-\boldsymbol{z}^*\right\|_{\nu}= \frac{\sqrt{2}\sigma M_{\boldsymbol{\mu}}(\|A\|e_1+\epsilon\nu^{\frac{1}{2}})}{1-\rho(\delta,\sigma,e_1)}.
\end{alignat}
\end{theorem}
\begin{IEEEproof}
See Appendix \ref{appendix:2}.
\end{IEEEproof}

We now discuss the impacts of gradient error $e_1$. The first term on the RHS of (\ref{convergeRegion})
$$\frac{\sqrt{2}\sigma M_{\boldsymbol{\mu}}\|A\|e_1}{1-\rho(\delta,\sigma,e_1)}$$
is nearly proportional to $e_1$, the approximate gradient error. Moreover, $\rho(\delta,\sigma,e_1)$ is a monotonically decreasing function of $e_1$, which again indicates that (\ref{convergeRegion}) is strictly decreasing with $e_1$. Reducing $e_1$ through an improved gradient evaluation will result in a better solution (i.e., a smaller sub-optimality gap) returned by dynamics (\ref{MultiIter:rewrite}).

\section{Numerical Results}
\label{sec:numerical}

% In this section, we demonstrate our improvement over the previous method \cite{zhou2019accelerated} through numerical results.
% \subsection{Experiment Setup}

We consider two test systems. The first is the IEEE $123$-bus test feeder. The second is a $4,518$-node test feeder (where we adopt the convention that a ``node'' refers to a certain phase at a certain bus), which is composed of the primary side of the IEEE $8,500$-node system and the ERPI Ckt7 system (with the secondary side lumped into primary-side transformers). Both feeders have unbalanced three phases. For the proposed hierarchical implementation, the $123$-bus feeder is clustered as shown in Figure \ref{fig:ieee123}. The clustering of the $4,518$-node feeder follows \cite{zhou2019accelerated}. For simplicity, the loads on all unclustered buses are set as fixed.

We make the following modifications to the original network models on the IEEE PES website (https://cmte.ieee.org/pes-testfeeders/resources/):
\begin{enumerate}
    \item All the loads are treated as constant-power loads. The delta-connected loads are omitted, and detailed models of capacitors, regulators, and breakers are not simulated.
    \item We double the original load data (used as the nominal case) of the $123$-bus feeder to create scenarios with serious under-voltage issues.
\end{enumerate}
\begin{figure*}
\begin{minipage}[h]{0.49\linewidth}
\centering
\includegraphics[width=1.00\columnwidth]{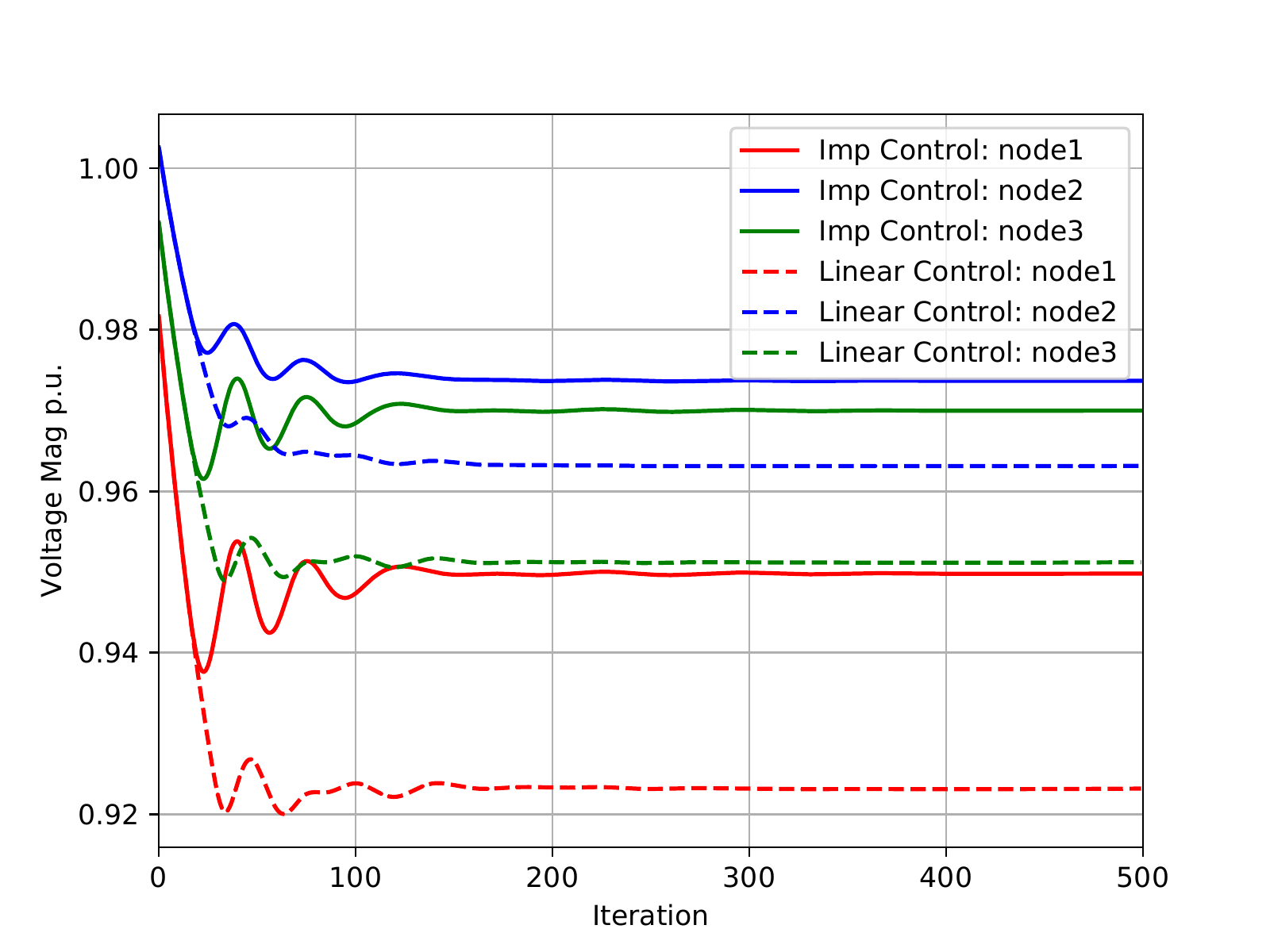}
% \caption{The voltage change at nodes of three different subtrees in the 4518-node network.}
% \label{fig:voltage_convergence_11000}
\end{minipage}
\begin{minipage}[h]{0.49\linewidth}
\centering
\includegraphics[width=1.00\columnwidth]{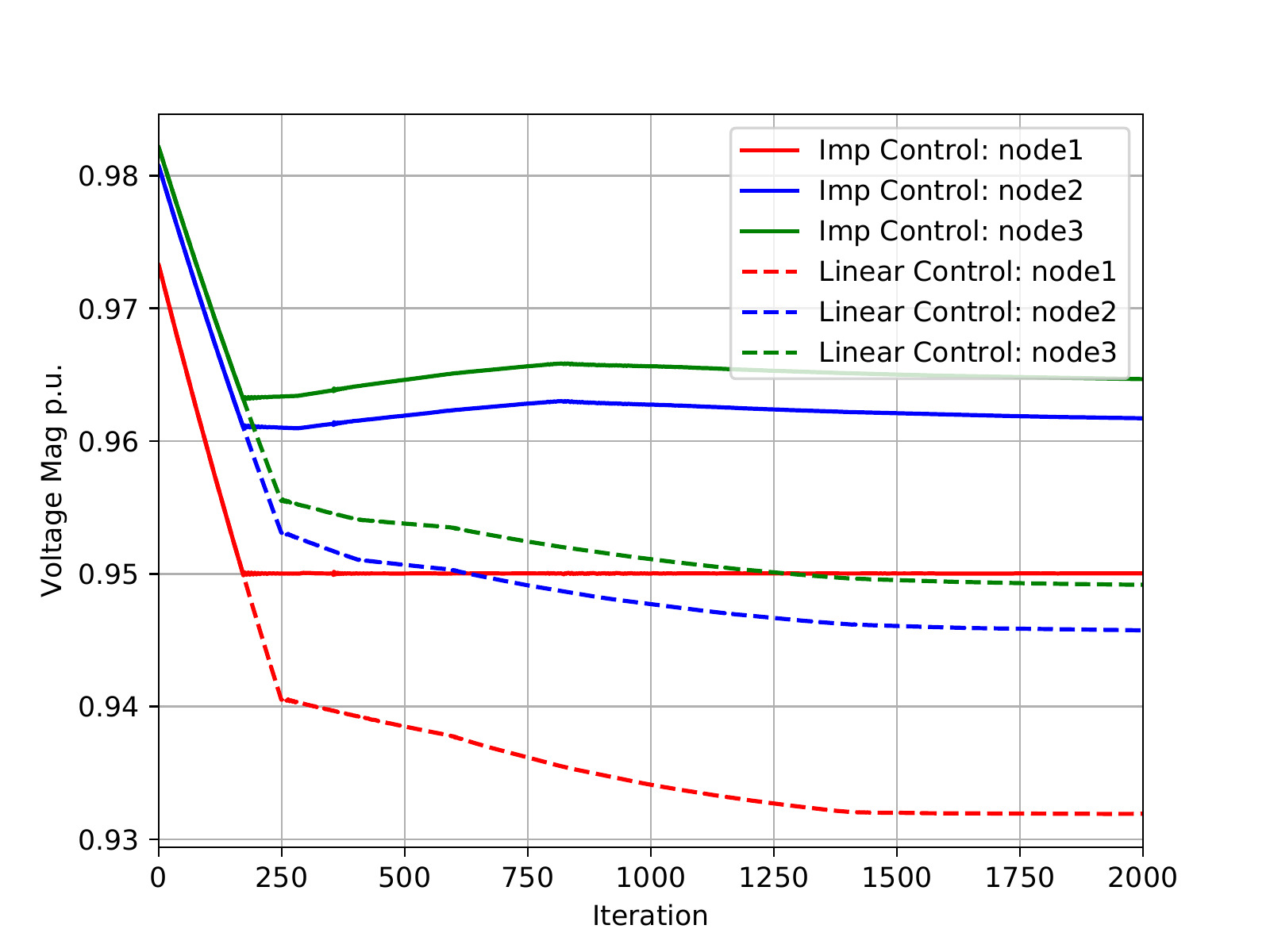}
\end{minipage}
\caption{The change of voltages at nodes from different subtrees in the IEEE $123$-bus feeder (left) and the $4,518$-node feeder (right).}
\label{fig:voltage_convergence}
\end{figure*} 

\begin{table*}[h]
    \centering
    \begin{tabular}{c|ccccc|cc}
\toprule
  \diagbox{Method}{Time (s)}{Subtrees} & \makecell{Subtree 1\\ Devices: 357} & \makecell{Subtree 2\\Devices: 222}& \makecell{Subtree 3\\Devices: 310} & \makecell{Subtree 4\\Devices: 154} & \makecell{All subtrees\\Devices: 1043} & Total Time\\
\midrule
 Improved control & $225.02$ & $141.02$ & $204.64$ & $53.08$ & $623.76$ & $1,003.78$ \\
\hline 
% \hline
 Linear control (hierarchical \cite{zhou2019accelerated}) & $198.22$ &  $103.47$ &  $168.52$ &  $52.96$ &  $523.17$ &  $570.49$ \\
\hline
 Linear control (centralized \cite{dall2016optimal}) &  -- &  --  &  --  &  -- & -- & $2,229.93$ \\
\bottomrule
\end{tabular}
    \caption{Run $2,000$ iterations on the $4,518$-node test feeder. ``Improved control'' refers to the hierarchical Algorithm \ref{algCCRC} with improved gradient evaluation (\ref{MultiImpGrad}), and ``linear control'' uses the lossless linear BFM (\ref{linearmodel}) for the primal-dual algorithm in a hierarchical or centralized manner. The number after ``devices'' refers to the number of controllable loads (buses).}
    \label{tab:running_time}
\end{table*}

For each phase $\varphi$ and bus $h$, we denote its nominal net power injection (which is negative for a load) by $(\underline{p}_{h}^{\varphi}, \underline{q}_{h}^{\varphi})$. The feasible power injection regions are defined as $\mathcal{Y}_{h}^{\varphi}=\left\{\left(p_{h}^{\varphi}, q_{h}^{\varphi}\right) \mid \underline{p}_{h}^{\varphi} \leqslant p_{h}^{\varphi} \leqslant 0.3\underline{p}_{h}^{\varphi} <0, \underline{q}_{h}^{\varphi} \leqslant q_{h}^{\varphi} \leqslant 0.3\underline{q}_{h}^{\varphi}<0\right\}$. The objective function in OPF problem (\ref{MultiOPF}) is defined as $ f_{h}^{\varphi}\left(p_{h}^{\varphi}, q_{h}^{\varphi}\right)=(p_h^{\varphi}-\underline{p}_h^{\varphi})^2+(q_h^{\varphi}-\underline{q}_h^{\varphi})^2$ for each controllable load, which aims to minimize the disutility caused by the deviation from the nominal loads.

For each test feeder, we fix the voltage magnitude of the root (slack) bus at $1.05$ per unit (p.u.), and set the lower and upper bounds for safe voltage at $0.95$ p.u. and $1.05$ p.u., respectively. For the $123$-bus feeder, the step sizes for primal and dual updates are chosen as $\sigma_{u}=1.5\times 10^{-2}$ and $\sigma_{\mu}=5.7\times 10^{-3}$, respectively. For the $4,518$-node feeder, we set $\sigma_{u}=7\times 10^{-4}$ and $\sigma_{\mu}=3\times 10^{-5}$. Both of them are chosen empirically such that $\rho(\delta,\sigma,e_1)<1$.

For each given power injection $\boldsymbol{u}:=[\boldsymbol{p}^{\top},\boldsymbol{q}^{\top}]^{\top}$, the OpenDSS platform is used to simulate the three-phase unbalanced power flow. The time for solving power flow by OpenDSS and interacting with OpenDSS is counted in our computational efficiency analysis. The nodal voltages $\boldsymbol{V}(\boldsymbol{u})$, line power flows $\boldsymbol{\Lambda}$, and currents $\boldsymbol{I}$ can be read from the simulator. We then convert $(\boldsymbol{V}(\boldsymbol{u}),\boldsymbol{\Lambda},\boldsymbol{I})$ to $(\boldsymbol{v}(\boldsymbol{u}),\boldsymbol{S},\boldsymbol{\ell})$, which are the variables in the BFM (\ref{distflow}). In the implementation, if we observe $\overline{\mu}_{j}^{\phi}(t)-\underline{\mu}_{j}^{\phi}(t)=0$ at any time $t$, we skip the associated gradient computation as it has no impact on (\ref{eq:decomposed_dvdp}) and (\ref{eq:decomposed_dvdp_unclustered}). The Python 3.7 programs for the proposed OPF algorithm and the OpenDSS simulation are run on a laptop equipped with Intel Core i7-9750H CPU @ 2.6 GHz, 16 GB RAM, and Windows 10 Professional OS.

\subsection{Voltage Safety}

The voltage magnitudes at all the nodes of each of the two test feeders are plotted in Figure \ref{fig:voltage_scatter}. There are three cases demonstrated in the figure. The first case, referred to as ``no control", takes the nominal power injections $(\underline{\boldsymbol{p}},\underline{\boldsymbol{q}})$, in which severe under-voltage violation is observed. The second case, referred to as ``linear control", applies the primal-dual algorithm proposed in literature \cite{zhou2019accelerated} for voltage regulation based on simplified gradient evaluation using lossless linearized BFM (\ref{linearmodel}). It enhances voltage safety to some extent, but still leaves many nodes below the lower limit. The third case, referred to as ``improved control", implements the proposed Algorithm \ref{algCCRC} based on the improved gradient evaluation (\ref{MultiImpGrad}), which successfully lifts all the nodes into the safe voltage range. This result verifies our analysis in Section \ref{sec:multi:model} that our method can prevent the voltage violation caused by ignoring the power loss during the previous model linearization. 

\subsection{Computational Efficiency}

Figure \ref{fig:voltage_convergence} shows the change of voltages at nodes from three different subtrees of each of the IEEE $123$-bus feeder and the $4,518$-node feeder. The dashed lines, corresponding to ``linear control", converge to a point where under-voltage violation occurs. On the contrary, the solid lines, corresponding to ``improved control", verify again that our method can prevent such voltage violation. Moreover, Table \ref{tab:running_time} displays the time to run $2,000$ iterations (which is more than needed for convergence as shown in Figure \ref{fig:voltage_convergence}) in the $4,518$-node feeder. In the table, the proposed ``improved control'' is compared with both the hierarchical and centralized implementations of ``linear control'' using the lossless model (the implementations in literature \cite{zhou2019accelerated} and \cite{dall2016optimal}, respectively). We observe that the proposed hierarchical algorithm with improved gradient evaluation saves over $50\%$ time than the centralized method, even if the latter used a simpler lossless model. Compared to the previous hierarchical method using the lossless model, the proposed method with improved gradient evaluation is slower by $14\%$ to $36\%$ in different subtrees and by $76\%$ in the whole network, due to the existence of bottleneck subtrees in our sequential simulation environment. Such increase in computation time can be alleviated by a real parallel implementation over subtrees. It is also an acceptable compromise considering the improved voltage safety of the proposed method.

\subsection{Comparison with SDP relaxation}
\begin{table}[!t]
\centering
\begin{tabular}{c|c|c|c|c|c}
\toprule %\hline 
\multirow{2}{*}{\text {Network}} & \multicolumn{2}{c|}{\text {Imp control}} & \multicolumn{2}{c|}{\text {SDP \cite{gan2014convex,zhao2017optimal}}} &\multirow{2}{*}{\text {Speedup}} \\
\cline { 2 - 5 } & \text{time (s)} & \text {solve?} & \text{time (s)} & \text {solve?} \\
\hline 
IEEE 37-bus & 1.20  & Y & 2.36  & Y & $\times$1.97 \\
IEEE 37-bus* & 1.81  & Y & 4.70  & N  & $\times$2.60 \\
IEEE 123-bus & 3.23 & Y & 35.25 & Y & $\times$10.91  \\
IEEE 123-bus* & 4.87 & Y & 59.38  & N & $\times$12.19  \\
\bottomrule
\end{tabular}
\caption{Comparison of ``improved control'' with SDP relaxation (* means  constraints \eqref{MultiOPF:v} are binding). ``Y'' means the OPF is solved, while ``N'' means not solved.}
\label{tab:compareSDP}
\end{table}

We compare the proposed ``improved control'' by running $300$ iterations (which is more than needed for convergence as shown in the left panel of Figure \ref{fig:voltage_convergence}), and compare its result with SDP relaxation in IEEE $37$-bus feeder and IEEE $123$-bus feeder with the same setup. We omit $4,518$-node feeder due to the extremely high computational burden of SDP relaxation. The generic optimization solver \emph{sedumi} \cite{jfsturn1999sedumi} based on Matlab 2016a is used to solve the SDP relaxation of OPF problem (\ref{MultiOPF}). Detailed models about SDP relaxation of (\ref{MultiOPF}) are referred to \cite{gan2014convex,zhao2017optimal}. Table \ref{tab:compareSDP} displays the results for cases where voltage constraints (\ref{MultiOPF:v}) are not binding or binding, by setting different nominal load data. 
In the case that voltage constraints are not binding, both methods can solve OPF. The proposed ``improved control'' achieves $\times1.97$ speedup in IEEE $37$-bus feeder and $\times10.91$ speedup in IEEE $123$-bus feeder. However, in the case that voltage constraints are binding, SDP relaxation fails to solve the OPF (\ref{MultiOPF}) due to violation of the rank-one constraint (\ref{disflow::rank}) \cite{gan2015exact}. The proposed improved control method still works in this case by returning a feasible and near-optimal  solution. Moreover, the proposed method achieves $\times2.60$ speedup in IEEE $37$-bus feeder and $\times12.19$ speedup in IEEE $123$-bus feeder.

\section{Conclusion}
\label{sec:conclusion}
We proposed a more accurate gradient evaluation method for three-phase unbalanced distribution networks. The resultant gradients with a blocked structure allow us to further design a scalable  hierarchical OPF algorithm. The proposed method prevents the previous voltage violation induced by the negligence of power loss in a linearized model, while achieving satisfactory computational efficiency, as verified by our theoretical analyses and numerical results. 

In the future, we shall further improve the accuracy of gradient evaluation with more efficient computation. We will devote more effort to theoretical analyses of convergence and optimality of the proposed algorithm. Besides characterizing algorithm performance with respect to a saddle point or a KKT point (which is often locally optimal for a nonconvex problem), we are interested in exploring its global sub-optimality bound. We are also adapting the proposed method to a time-varying setting that incorporates time-coupling constraints (e.g., for energy storage) to enable fast online solutions.

{\appendices
\section{Proof of Lemma \ref{lemma::Vsafe:three}}
\label{appendix::1}
Take the difference between the intermediate model (\ref{disflow::v}), (\ref{BTdisflow}) and the linearized lossless model (\ref{linearmodel}):
\begin{subequations}\label{MultiError}
\begin{alignat}{2} 
\hat{v}_{j}-v_{j}&= \hat{v}_{i}^{\Phi_{ij}}-v_{i}^{\Phi_{ij}}-\Big(\left(\hat{S}_{i j}-S_{ij}\right) z_{i j}^{H} \nonumber \\
& \quad +z_{i j}\left( \hat{S}_{i j}-S_{ij}\right)^{H}\Big)- z_{i j} \ell_{i j} z_{i j}^{H}, \label{MultiError::v}\\
\hat{S}_{ij}-S_{ij}&=\gamma^{\Phi_{ij}}\operatorname{Diag}(\hat{\Lambda}_{ij}-\Lambda_{ij}), \label{MultiError::SL}\\
\hat{\Lambda}_{ij}-\Lambda_{ij}&=\sum_{k:(j,k)\in \mathcal{E}}(\hat{\Lambda}_{jk}^{\Phi_j}-\Lambda_{jk}^{\Phi_j})-\operatorname{diag}(z_{ij}\ell_{ij}). \label{MultiError::PQ}
\end{alignat}
\end{subequations}

Notice that there exists a power flow error term $\hat{S}_{ij}-S_{ij}$ on the RHS of (\ref{MultiError::v}), and this error is calculated in (\ref{MultiError::SL})--(\ref{MultiError::PQ}). Following the tree structure of the network, we recursively calculate $\hat{S}_{ij}-S_{ij}$ through (\ref{MultiError::SL})--(\ref{MultiError::PQ}) until reaching leaves: 
\begin{subequations} \label{Serror}
\begin{alignat}{2} 
\hat{S}_{ij}-S_{ij}=&\sum_{k:(j,k)\in \mathcal{E}}(\hat{S}_{jk}^{\Phi_j}-S_{jk}^{\Phi_j}) \nonumber\\
&-\gamma^{\Phi_{ij}}\operatorname{Diag}(\operatorname{diag}(z_{ij}\ell_{ij})),\label{Serror:ij} \\
\hat{S}_{jk}-S_{jk}=&\sum_{l:(k,l)\in \mathcal{E}}(\hat{S}_{kl}^{\Phi_k}-S_{kl}^{\Phi_k}) \nonumber\\ 
&-\gamma^{\Phi_{jk}}\operatorname{Diag}(\operatorname{diag}(z_{jk}\ell_{jk})),\label{Serror:jk} \\
% \hat{S}_{kt}-S_{kt}&=\sum_{h:(t,h)\in \mathcal{E}}(\hat{S}_{th}-S_{th})-\gamma\operatorname{Diag}(\operatorname{diag}(z_{kt}\ell_{kt})),\label{Serror:kt}\\
\qquad \dots \nonumber \\
\hat{S}_{th}-S_{th}=&-\gamma^{\Phi_{th}}\operatorname{Diag}(\operatorname{diag}(z_{th}\ell_{th})), \text{ for leaf buses $h$}.\label{Serror:th}
\end{alignat}
\end{subequations}

\begin{figure}
\centering
\includegraphics[width=0.90\columnwidth]{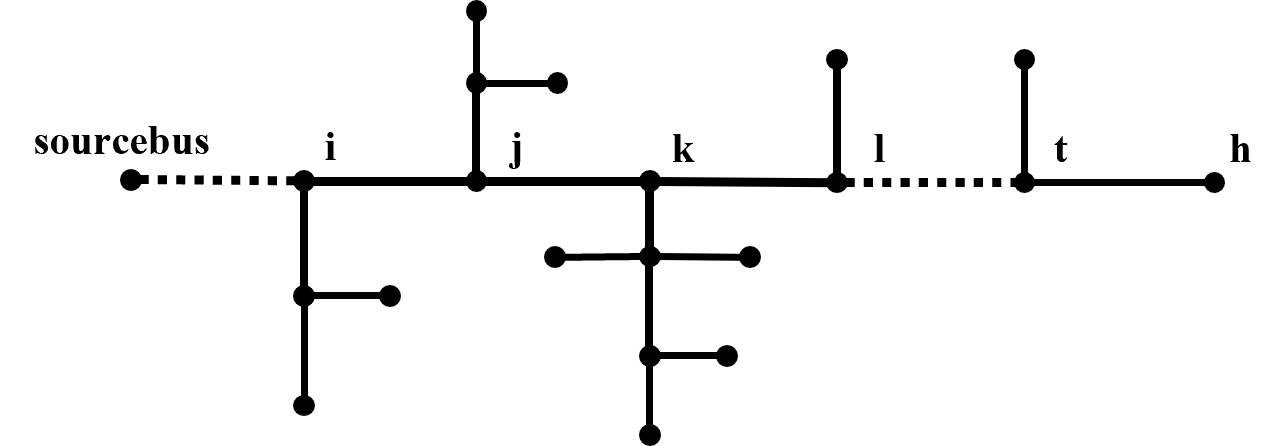}
\caption{Index the buses in the proof of \textit{Lemma \ref{lemma::Vsafe:three}} (Appendix).} \label{fig:treegraph}
\end{figure}
The bus indices in (\ref{Serror}) are illustrated in Figure \ref{fig:treegraph}. Recursively substituting such equations for each layer of the network into its upper layer, i.e., (\ref{Serror:th}), $\dots$, into (\ref{Serror:jk}), then into (\ref{Serror:ij}), we have:
\begin{alignat}{2} \label{Serror:ij:2}
& \hat{S}_{ij} -S_{ij}= \nonumber\\
& \quad-\sum_{k:(j,k)\in \mathcal{E}}\sum_{l:(k,l)\in \mathcal{E}}\dots\sum_{h:(t,h)\in \mathcal{E}}\gamma^{\Phi_{th}}\operatorname{Diag}(\operatorname{diag}(z_{th}\ell_{th}))\nonumber\\
& \quad-...-\sum_{k:(j,k)\in \mathcal{E}}\sum_{l:(k,l)\in \mathcal{E}}\gamma^{\Phi_{kl}}\operatorname{Diag}(\operatorname{diag}(z_{kl}\ell_{kl})) \nonumber\\
& \quad-\sum_{k:(j,k)\in \mathcal{E}}\gamma^{\Phi_{jk}}\operatorname{Diag}(\operatorname{diag}(z_{jk}\ell_{jk})) \nonumber \\ 
& \quad-\gamma^{\Phi_{ij}}\operatorname{Diag}(\operatorname{diag}(z_{ij}\ell_{ij})).
\end{alignat}

In fact, the RHS of (\ref{Serror:ij:2}) is related to the power losses downstream of bus $i$. Define $M_{ij}:=S_{ij}-\hat{S}_{ij}$, which can be written as:
\begin{alignat}{2}
M_{ij} & =\gamma^{\Phi_{ij}}\operatorname{Diag}(\operatorname{diag}(z_{ij}\ell_{ij}))  \nonumber\\
& +\sum_{(\alpha,\beta)\in\operatorname{down}(j)}\gamma^{\Phi_{\alpha \beta}}\operatorname{Diag}(\operatorname{diag}(z_{\alpha\beta}\ell_{\alpha\beta})).
\end{alignat}
From the tree structure, we have:
\begin{alignat}{2} \label{MijRewrite}
M_{ij}=\sum_{k:(j,k)\in \mathcal{E}}M_{jk}+\gamma^{\Phi_{ij}}\operatorname{Diag}(\operatorname{diag}(z_{ij}\ell_{ij})).
\end{alignat}

The relationship (\ref{MijRewrite}) indicates a layered relationship in the tree network. A backward sweep from the leaves to the root will return $M_{ij}, ~\forall (i,j)\in \mathcal{E}$. Based on $M_{ij}$, the voltage error (\ref{MultiError::v}) can be written as:
\begin{alignat}{2} \label{VError::ij}
\hat{v}_{j}-v_{j}= \hat{v}_{i}^{\Phi_{ij}}-v_{i}^{\Phi_{ij}}+(M_{ij}z_{ij}^{H}+z_{ij}M_{ij}^{H})-z_{ij}\ell_{ij}z_{ij}^{H}.
\end{alignat}

Equation (\ref{VError::ij}) implies an accumulated voltage error from the root to leaf buses. A forward sweep from the root to leaves will return the voltage errors at all buses $h\in\mathcal{N}$ as:
\begin{alignat}{2} \label{VError::h}
\hat{v}_{h}-v_{h}=\sum_{(\zeta,\xi)\in \mathbb{P}_{h}}\left(M_{\zeta\xi}z_{\zeta\xi}^{H}+z_{\zeta\xi}M_{\zeta\xi}^{H}-z_{\zeta\xi} \ell_{\zeta\xi}z_{\zeta\xi}^{H}\right).
\end{alignat}

This completes the proof of \textit{Lemma 1}.

\section{Proof of Theorem \ref{theom:1}}
\label{appendix:2}

The KKT condition (\ref{KKT}) implies:
\begin{subequations} \label{stationarypoint}
\begin{alignat}{2}
&\boldsymbol{u}^*=\bigg[\boldsymbol{u}^*-\sigma\bigg(\nabla f(\boldsymbol{u}^*)+\frac{\partial \boldsymbol{v}(\boldsymbol{u}^*)}{\partial \boldsymbol{u}}^{\top}A^{\top}\boldsymbol{\mu}^*\bigg)\bigg]_{\mathcal{Y}}, \\
&\boldsymbol{\mu}^*=\bigg[\boldsymbol{\mu}^*+\nu\sigma\left(A\boldsymbol{v}(\boldsymbol{u}^*)-d\right)\bigg]_+.
\end{alignat}
\end{subequations}

Comparing (\ref{MultiIter:rewrite}) and (\ref{stationarypoint}), and by the non-expansion of projection and the definition (\ref{def:Delta}) of $\Delta$, we have:
\begin{subequations}
\begin{alignat}{2}
&\left\|\boldsymbol{u}(t+1)-\boldsymbol{u}^*\right\|^2\leq \nonumber\\
&\left\|(1-\sigma\Delta)(\boldsymbol{u}(t)-\boldsymbol{u}^*)-\sigma \frac{\partial \widetilde{\boldsymbol{v}}(\boldsymbol{u}(t))}{\partial \boldsymbol{u}}^{\top}A^{\top}\left(\boldsymbol{\mu}(t)-\boldsymbol{\mu}^*\right) \right. \nonumber\\
&\qquad \left. -\sigma \left(\frac{\partial \widetilde{\boldsymbol{v}}(\boldsymbol{u}(t))}{\partial \boldsymbol{u}}-\frac{\partial \boldsymbol{v}(\boldsymbol{u}(t))}{\partial \boldsymbol{u}}\right)^{\top}A^{\top}\boldsymbol{\mu}^*\right\|^2,  \\
&\left\|\boldsymbol{\mu}(t+1)-\boldsymbol{\mu}^*\right\|^2 \leq \left\|(1-\nu\sigma\epsilon)\left(\boldsymbol{\mu}(t)-\boldsymbol{\mu}^*\right)-\nu\sigma\epsilon\boldsymbol{\mu}^*\right. \nonumber\\
& \qquad\left.+\nu\sigma A \left(\boldsymbol{v}(\boldsymbol{u}(t))-\boldsymbol{v}(\boldsymbol{u}^*)\right)\right\|^2.
\end{alignat}
\end{subequations}

Moreover, the definitions (\ref{def::nonlinear}) and Assumption \ref{ass:boundgradient} imply:
\begin{alignat}{2} \label{bounds::1}
&\left\|\frac{\partial \widetilde{\boldsymbol{v}}(\boldsymbol{u}(t))}{\partial \boldsymbol{u}}^{\top}\left(\boldsymbol{\mu}(t)-\boldsymbol{\mu}^*\right)\right\|\leq (L_{\boldsymbol{v}}(\delta)+e_1)\left\|\boldsymbol{\mu}(t)-\boldsymbol{\mu}^*\right\|, \nonumber\\
&\qquad \left\|\boldsymbol{v}(\boldsymbol{u}(t))-\boldsymbol{v}(\boldsymbol{u}^*)\right\| \leq L_{\boldsymbol{v}}(\delta)\left\|\boldsymbol{u}(t)-\boldsymbol{u}^*\right\|, \nonumber\\
&\qquad\left\|\left(\frac{\partial \widetilde{\boldsymbol{v}}(\boldsymbol{u}(t))}{\partial \boldsymbol{u}}-\frac{\partial \boldsymbol{v}(\boldsymbol{u}(t))}{\partial \boldsymbol{u}}\right)^{\top}\boldsymbol{\mu}^*\right\|\leq e_1M_{\boldsymbol{\mu}}
\end{alignat}
and
\begin{alignat}{2} \label{bounds::2}
&\left\|\boldsymbol{v}(\boldsymbol{u}(t))-\boldsymbol{v}(\boldsymbol{u}^*)-\frac{\partial \widetilde{\boldsymbol{v}}(\boldsymbol{u}(t))}{\partial \boldsymbol{u}}^{\top}\left(\boldsymbol{u}(t)-\boldsymbol{u}^*\right)\right\|  \nonumber\\
& \qquad \leq(M_{\boldsymbol{v}}(\delta)\delta+e_1)\left\|\boldsymbol{u}(t)-\boldsymbol{u}^*\right\|.
\end{alignat}

By using the norm $\left\|\boldsymbol{z}\right\|_{\nu}^2$ and combining the bounds in (\ref{bounds::1}) and (\ref{bounds::2}), we have:
\begin{alignat}{2}\label{projection:z}
&\quad\left\|\boldsymbol{z}(t+1)-\boldsymbol{z}^*\right\|_{\nu}^2 \leq \nonumber\\
&\left((1-\sigma\Delta)^2+\nu\sigma^2\|A\|^2L_{\boldsymbol{v}}^2(\delta)\right)\left\|\boldsymbol{u}(t)-\boldsymbol{u}^*\right\|^2 \nonumber\\
&+\left((1-\nu\sigma\epsilon)^2+\nu\sigma^2\|A\|^2(L_{\boldsymbol{v}}(\delta)+e_1)^2\right)\nu^{-1}\left\|\boldsymbol{\mu}(t)-\boldsymbol{\mu}^*\right\|^2 \nonumber \\
& +2\sigma \big[(1-\nu\sigma\epsilon)(M_{\boldsymbol{v}}(\delta)\delta+e_1)+\sigma|\Delta-\nu\epsilon|(L_{\boldsymbol{v}}(\delta)+e_1)\big] \nonumber\\
&\cdot\|A\|\left\|\boldsymbol{u}(t)-\boldsymbol{u}^*\right\|\left\|\boldsymbol{\mu}(t)-\boldsymbol{\mu}^*\right\| \nonumber \\
&+2\sigma M_{\boldsymbol{\mu}}\epsilon \nonumber\\
&\cdot\big[\nu\sigma\|A\|L_{\boldsymbol{v}}(\delta)\left\|\boldsymbol{u}(t)-\boldsymbol{u}^*\right\|+(1-\nu\sigma\epsilon)\left\|\boldsymbol{\mu}(t)-\boldsymbol{\mu}^*\right\|\big] \nonumber\\
&+2\sigma M_{\boldsymbol{\mu}}e_1 \|A\| \nonumber\\
&\cdot\big[(1-\sigma\Delta)\left\|\boldsymbol{u}(t)-\boldsymbol{u}^*\right\|+\sigma\|A\|(L_{\boldsymbol{v}}(\delta)+e_1)\left\|\boldsymbol{\mu}(t)-\boldsymbol{\mu}^*\right\|\big] \nonumber \\
&+\sigma^2 M_{\boldsymbol{\mu}}^2(\|A\|^2e_1^2+\epsilon^2\nu).
\end{alignat}

Furthermore, by using Young's inequality, we have:
\begin{alignat}{2} \label{crossterm:1}
&\quad\left\|\boldsymbol{u}(t)-\boldsymbol{u}^*\right\|\left\|\boldsymbol{\mu}(t)-\boldsymbol{\mu}^*\right\| \leq \nonumber\\
&\frac{1}{2}\sqrt{\nu}\left(\left\|\boldsymbol{u}(t)-\boldsymbol{u}^*\right\|^2+\nu^{-1}\left\|\boldsymbol{\mu}(t)-\boldsymbol{\mu}^*\right\|^2\right).
\end{alignat}
And 
\begin{subequations} \label{crossterm:2}
\begin{alignat}{2}
&2\sqrt{2}\sigma M_{\boldsymbol{\mu}}\epsilon\nu^{\frac{1}{2}}\bigg[\frac{\nu^{\frac{1}{2}}\sigma\|A\|L_{\boldsymbol{v}}(\delta)}{\sqrt{2}}\left\|\boldsymbol{u}(t)-\boldsymbol{u}^*\right\|  \nonumber\\
&\quad+\frac{(1-\nu\sigma\epsilon)}{\sqrt{2}}\nu^{-\frac{1}{2}}\left\|\boldsymbol{\mu}(t)-\boldsymbol{\mu}^*\right\|\bigg] \nonumber\\
& \leq 2\sqrt{2}\sigma M_{\boldsymbol{\mu}}\epsilon\nu^{\frac{1}{2}}\big[\nu\sigma^2\|A\|^2L_{\boldsymbol{v}}^2(\delta)\left\|\boldsymbol{u}(t)-\boldsymbol{u}^*\right\|^2 \nonumber\\
&\quad+(1-\nu\sigma\epsilon)^2\nu^{-1}\left\|\boldsymbol{\mu}(t)-\boldsymbol{\mu}^*\right\|^2\big]^{\frac{1}{2}} ,\\
&2\sqrt{2}\sigma M_{\boldsymbol{\mu}}e_1\|A\|\bigg[\frac{(1-\sigma\Delta)}{\sqrt{2}}\left\|\boldsymbol{u}(t)-\boldsymbol{u}^*\right\| \nonumber\\
&\quad+\frac{\nu^{\frac{1}{2}}\sigma\|A\|(L_{\boldsymbol{v}}(\delta)+e_1)}{\sqrt{2}}\nu^{-\frac{1}{2}}\left\|\boldsymbol{\mu}(t)-\boldsymbol{\mu}^*\right\|\bigg] \nonumber\\
&\leq 2\sqrt{2}\sigma M_{\boldsymbol{\mu}}e_1\|A\|\big[(1-\sigma\Delta)^2\left\|\boldsymbol{u}(t)-\boldsymbol{u}^*\right\|^2 \nonumber\\
&\quad+\nu\sigma^2\|A\|^2(L_{\boldsymbol{v}}(\delta)+e_1)^2\nu^{-1}\left\|\boldsymbol{\mu}(t)-\boldsymbol{\mu}^*\right\|^2\big]^{\frac{1}{2}}.
\end{alignat}
\end{subequations}

Then by combining (\ref{crossterm:1}), (\ref{crossterm:2}) with (\ref{projection:z}), and by the definition $\rho(\delta,\sigma,e_1)$ in (\ref{convergeSpeed}), we can show that:
\begin{alignat}{2}
&\left\|\boldsymbol{z}(t+1)-\boldsymbol{z}^*\right\|_{\nu}^2 \leq \rho^2(\delta,\sigma,e_1)\left\|\boldsymbol{z}(t)-\boldsymbol{z}^*\right\|_{\nu}^2 \nonumber\\
&+2\sqrt{2}\sigma M_{\boldsymbol{\mu}}(\epsilon\nu^{\frac{1}{2}}+e_1 \|A\|)\rho(\delta,\sigma,e_1)\left\|\boldsymbol{z}(t)-\boldsymbol{z}^*\right\|_{\nu} \nonumber\\
&\quad + \left(\sqrt{2}\sigma M_{\boldsymbol{\mu}}(\|A\|e_1+\epsilon\nu^\frac{1}{2})\right)^2.
\end{alignat}

It is equivalent to:
\begin{alignat}{2}
\left\|\boldsymbol{z}(t+1)-\boldsymbol{z}^*\right\|_{\nu} \leq \rho(\delta,\sigma,e_1)\left\|\boldsymbol{z}(t)-\boldsymbol{z}^*\right\|_{\nu} \nonumber\\
+\sqrt{2}\sigma M_{\boldsymbol{\mu}}(\|A\|e_1+\epsilon\nu^{\frac{1}{2}}) \nonumber\\
\leq\rho^{t}(\delta,\sigma,e_1)\left\|\boldsymbol{z}(1)-\boldsymbol{z}^*\right\|_{\nu}\nonumber\\
+\frac{1-\rho^{t}(\delta,\sigma,e_1)}{1-\rho(\delta,\sigma,e_1)}\sqrt{2}\sigma M_{\boldsymbol{\mu}}(\|A\|e_1+\epsilon\nu^{\frac{1}{2}}).
\end{alignat}

Since $\rho(\delta,\sigma,e_1)<1$, it follows that:
\begin{alignat}{2}
\limsup_{t \rightarrow \infty}\left\|\boldsymbol{z}(t)-\boldsymbol{z}^*\right\|_{\nu}=\frac{\sqrt{2}\sigma M_{\boldsymbol{\mu}}(\|A\|e_1+\epsilon\nu^{\frac{1}{2}})}{1-\rho(\delta,\sigma,e_1)}. \nonumber
\end{alignat}

This completes the proof of \textit{Theorem 1}.
}

\vfill

\end{document}